\documentclass{article}

    \usepackage[final]{neurips_2022}

\usepackage[utf8]{inputenc} 
\usepackage[T1]{fontenc}    
\usepackage{hyperref}       
\usepackage{url}            
\usepackage{booktabs}       
\usepackage{amsfonts}       
\usepackage{nicefrac}       
\usepackage{microtype}      
\usepackage{xcolor}         
\usepackage{booktabs}

\usepackage{subfigure}

\usepackage{algorithmic}
\usepackage{graphicx}
\usepackage{color,colortbl}

\usepackage{amsmath}
\usepackage{amssymb}
\usepackage{mathtools}
\usepackage{amsthm}
\usepackage{multirow}
\usepackage{multicol}
\usepackage[capitalize,noabbrev]{cleveref}

\usepackage{comment}
\usepackage{optidef}
\usepackage{enumitem}

\usepackage{algorithm}

\usepackage{empheq}

\newcommand{\jd}[1]{\textcolor{blue}{#1}} 

\newcommand{\vect}[1]{\mathbf{#1}} 
\newcommand{\vects}[1]{\boldsymbol{#1}} 

\newcommand{\EVA}{NN-E}
\newcommand{\PSA}{NN-P}

\newcommand{\method}{\texttt{Neur2SP}}

\newcommand{\twostocprog}{2SP}
\newcommand{\extform}{\text{EF}}

\newcommand{\NNinput}{\ensuremath{\vects{\alpha}}}
\newcommand{\NNoutput}{\ensuremath{\beta}}

\title{\method: Neural Two-Stage Stochastic Programming}

\author{
	Justin Dumouchelle\thanks{These authors contributed equally.}\qquad
	Rahul Patel$^*$ \qquad
	Elias B. Khalil\thanks{Corresponding author: \href{mailto:khalil@mie.utoronto.ca}{\url{khalil@mie.utoronto.ca}}.} \qquad
	Merve Bodur \\
    Department of Mechanical \& Industrial Engineering, University of Toronto
}

\hypersetup{
pdftitle={Neur2SP: Neural Two-Stage Stochastic Programming},
pdfsubject={cs.LG, math.OC},
pdfauthor={Justin Dumouchelle, Rahul Patel, Elias B. Khalil, Merve Bodur},
pdfkeywords={Two-stage Stochastic Programming, Supervised Learning, Deep Learning},
}

\begin{document}
\maketitle

\begin{abstract}
Stochastic Programming is a powerful modeling framework for decision-making under uncertainty.
In this work, we tackle two-stage stochastic programs (\twostocprog{}s), the most widely used class of stochastic programming models.
Solving \twostocprog{}s exactly requires optimizing over an expected value function that is computationally intractable. 
Having a mixed-integer linear program (MIP) 
or a nonlinear program (NLP) in the second stage further aggravates the intractability, even when specialized algorithms that exploit problem structure are employed.
Finding high-quality (first-stage) solutions -- without leveraging problem structure -- can be crucial in such settings.
We develop~\method{}, a new method that approximates the expected value function via a neural network to obtain a surrogate model that can be solved more efficiently than the traditional extensive formulation approach. 
\method{} makes no assumptions about the problem structure, in particular about the second-stage problem, and can be implemented using an off-the-shelf MIP solver.
Our extensive computational experiments on four benchmark \twostocprog{} problem classes with different structures (containing MIP and NLP second-stage problems) demonstrate the efficiency (time) and efficacy (solution quality) of~\method.
In under 1.66 seconds,~\method{} finds high-quality solutions across all problems even as the number of scenarios increases, an ideal property that is difficult to have for traditional \twostocprog{} solution techniques. Namely, the most generic baseline method typically requires minutes to hours to find solutions of comparable quality.

\end{abstract}


\begin{figure}[h]
    \centering
    \includegraphics[width=\linewidth]{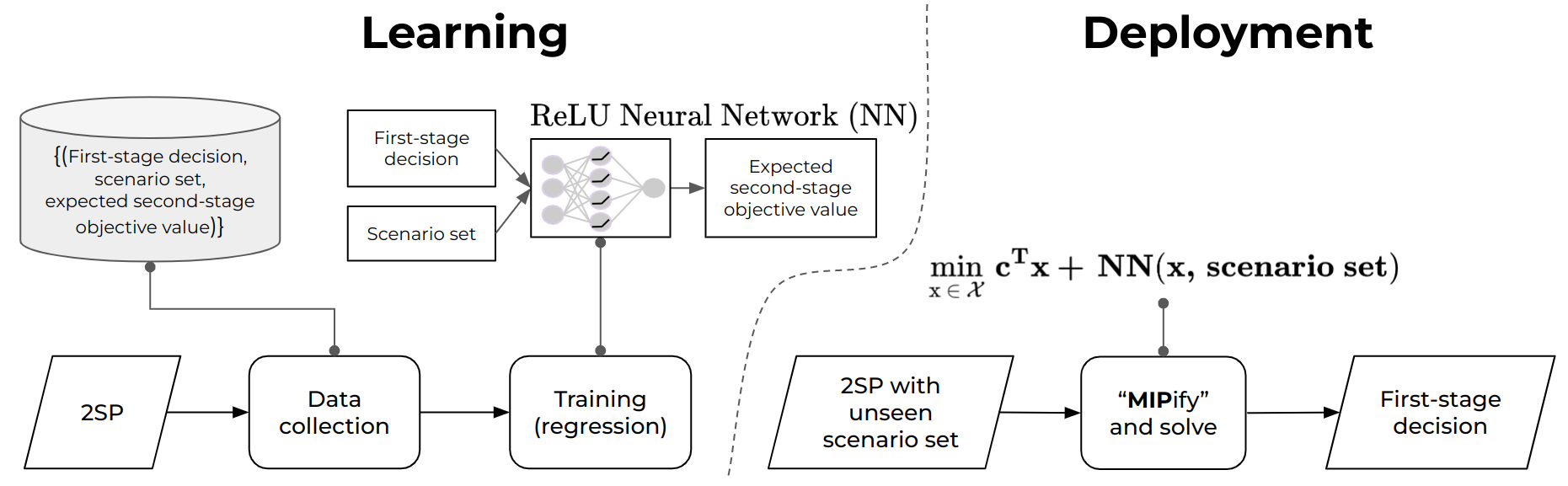}
    \caption{
    {
    Overview of~\method{}. The leftmost block is the input, namely, a \twostocprog{}. From the \twostocprog{}, we follow the data generation procedure from Section~\ref{sec:methodology_dg} to obtain a dataset consisting of tuples of (first-stage decision,  scenario set, corresponding expected second-stage objective value). We then train one of the learning models presented in Section~\ref{sec:methodology_lm} to predict the expected cost given a first-stage decision and scenario set.  The trained model is then embedded into a MIP using the procedure in Section~\ref{sec:methodology_embed} to obtain an approximate MIP (the ``MIPify'' step). Lastly, the approximate MIP is solved with an off-the-shelf MIP solver to obtain a first-stage 2SP solution.}
    }
    \label{fig:block-diagram}
\end{figure}

\section{Introduction}

Mathematical programming consists of a gamut of tools to solve optimization problems.
Under perfect information, i.e., when all the data is deterministic and known, many of these problems can be solved as a linear program (LP) or mixed-integer linear program (MIP).
However, in many cases, there is a need to deal with problems with partial information.
Stochastic programming is one such framework that allows us to incorporate uncertainty into decision-making.

In this work, we focus our attention on Two-stage Stochastic Programs (\twostocprog s).  
A \twostocprog{} involves two sets of decisions, namely the first-stage and second-stage (recourse) decisions, to be made before and after the uncertainty is realized, respectively. Given the (joint) probability distribution of the random parameters of the problem, the most common objective of \twostocprog{} is to optimize the expected value of the decisions.  
For example, in a two-stage stochastic facility location problem, first-stage decisions consist of which facilities should be built whereas second-stage decisions involve assigning customers to open facilities to meet their stochastic demand, and the overall objective is to minimize the sum of the cost of the first-stage decisions and the expected cost of the second-stage decisions.

\twostocprog{s} are usually solved via Sample Average Approximation (SAA), which limits the future uncertainty to a finite set of possible realizations (scenarios). The SAA approximation of a 2SP is a reduction to an equivalent deterministic problem and can be solved by the so-called \textit{extensive form}: a monolithic formulation where scenario copies of the  second-stage decision variables are created and linked to the first-stage decisions. However, even for small 2SPs, solving the extensive form may be intractable as it requires introducing a large number of (possibly integer) variables and (possibly nonlinear) constraints. 
As such, specialized algorithms are required.
If the second-stage problem assumes the form of an LP, then algorithms such as Benders' decomposition (also known as the L-shaped method) can be leveraged to efficiently solve the problem to optimality.
Unfortunately, in many practical applications of 2SP, the second-stage problem assumes the form of a MIP, for which specialized decomposition algorithms might not be efficient. The existence of continuous first-stage variables linked to the second-stage problem significantly increases the difficulty of solving such problems. This is exacerbated when the second-stage problem is nonlinear, for which no general and structure-agnostic solution strategy exists.

In this work, we propose~\method{}, a framework for constructing an {easier-to-solve} surrogate optimization problem for 2SP with the use of supervised deep learning. In a nutshell, a Rectified Linear Unit (ReLU) neural network is trained to approximate the second-stage objective value for a set of scenarios. Using MIP-representable activation functions such as the ReLU, the forward computation of the trained network can be embedded into a MIP.
The surrogate problem is then confined to optimizing \textit{only} first-stage decisions with respect to the first-stage objective function and the neural network approximation of the second-stage objective~\footnote{For a fixed first-stage solution obtained via this surrogate, an optimal second-stage decision can be obtained relatively quickly for each scenario if desired.}. Assuming a small and accurate neural network can be used, the surrogate problem is much smaller than the extensive form, and thus faster to solve.   
The entire procedure is summarized in \Cref{fig:block-diagram}.
Our main contributions are as follows: 
\begin{enumerate}
    \item \textbf{Novelty:} \method{} is the first generic machine learning approach for deriving a heuristic solution for \twostocprog{}. We introduce a highly parallelizable data collection procedure and show two separate neural models 
    which can be used to formulate a
    deterministic mixed-integer surrogate problem for 2SP;
    \item \textbf{Generality:}  \method{} can be used out-of-the-box for \twostocprog{s} with linear and nonlinear objectives and constraints as well as mixed-integer variables in both the first and second stages, all without using any problem structure, i.e., in a purely data-driven way;
    \item \textbf{Performance:} \method{} is shown to produce high-quality solutions significantly faster than the solely applicable general baseline method, the extensive form approach, for 
    a variety of benchmark problems, namely, stochastic facility location problem, an investment problem, a server location problem, and a pooling problem from chemical engineering.

\end{enumerate}

\section{Preliminaries}
\label{sec:preliminaries}
We introduce the \twostocprog{} setting and describe the MIP formulation for a ReLU activation function which is central to the surrogate model we propose in this work. Appendix~\ref{app:symbols} summarizes the notation used.


\subsection{Two-stage Stochastic Programming}
\label{sec:preliminaries_twosp}
A 2SP can be generally expressed as
$\min_{\vect x}\{\vect c^\intercal \vect x + \mathbb{E}_{\vect \xi}[Q(\vect x, \vects \xi)] : \vect x \in \mathcal X\},$
where $\vect c \in \mathbb{R}^n$ is the first-stage cost vector, $\vect x \in \mathbb{R}^n$ represents the first-stage decisions, $\mathcal X$ is the first-stage feasible set, and $\vects \xi$ is the vector of random parameters that follow a probability distribution $\mathbb{P}$ with support $\Xi$. 
The \emph{value function} $Q: \mathcal X \times \Xi \rightarrow \mathbb{R}$ returns the cost of optimal second-stage (recourse) decisions under realization $\vects \xi$ given the first-stage decisions of $\vect x$.  In many cases, as the $Q(\vect x, \vects \xi)$ is obtained by solving a mathematical program, evaluating the \emph{expected value function} $\mathbb{E}_{\vects \xi}[Q(\vect x, \vects \xi)]$ is intractable.  

To provide a more tractable formulation, the extensive form (\extform{}) is used.
Using a set of $K$ scenarios, $\vects \xi_1,\ldots,\vects \xi_K$, sampled from the probability distribution $\mathbb{P}$, $\extform(\vects \xi_1,\ldots,\vects \xi_K) \equiv \min_{\vect x}\{\vect c^\intercal \vect x + \sum_{k=1}^{K} p_k Q(\vect x, \vects \xi_k) :  \vect x \in \mathcal{X}\},$
where $p_k$ is the probability of scenario $\vects \xi_k$ being realized.  
If $Q(\vect x, \vects \xi) = \min_{\vect y} \{ F(\vect y,\vects \xi) : \vect y \in \mathcal{Y}(\vect x,\vects \xi) \}$, then $\extform{(\vects \xi_1,\ldots,\vects \xi_K)}$ can be expressed as 
$$\min_{\vect x,\vect y} \Bigg\{ \vect c^\intercal \vect x + \sum_{k=1}^{K} p_k F(\vect y_k,\vects \xi_k): \vect x \in \mathcal X, \vect y_k \in \mathcal{Y}(\vect x,\vects \xi_k) \, \forall k=1,\hdots, K \Bigg\},$$ which can be solved through standard deterministic optimization techniques.
However, the number of variables and  constraints of the \extform{} grows linearly with the number of scenarios.  
Furthermore, if $Q(\cdot, \cdot)$ is the optimal value of a MIP or an nonlinear program (NLP), the \extform{} model becomes significantly more challenging to solve as compared to the LP case, limiting its applicability even at small scale.



\subsection{Embedding Neural Networks into MIPs}
\label{sec:nn_to_mip}
Mathematically, an $\ell$-layer fully-connected neural network can be expressed as:
$\vect h^{1} = \sigma(W^{0} \NNinput + \vect b^{0})$; $\vect h^{m+1} = \sigma(W^{m}  \vect  h^{m} + \vect b^{m}), m = 1,\dots,\ell-1$; $\NNoutput = W^{\ell}  \vect h^{\ell} + \vect b^{\ell}$.
Here, $\NNinput \in\mathbb{R}^m$ is the input, $\NNoutput \in\mathbb{R}$ is the prediction, $\vect h^{i}\in\mathbb{R}^{d_i}$ is the $i$-th hidden layer, $W^{i}\in\mathbb{R}^{d_{i} \times d_{i+1}}$ is the matrix of weights from layer $i$ to $i+1$, $\vect b^{i}\in\mathbb{R}^{d_i}$ is the bias at the $i$-th layer, and $\sigma$ is a non-linear activation function, here the activation function is given by $\text{ReLU}(a) = \max\{0, a\}$ for $a\in\mathbb{R}$.

Central to~\method{} is the embedding of a trained neural network into a MIP.  Here, we present the formulation proposed by \citep{fischetti2018deep}. For a given hidden layer $m$, the $j$-th hidden unit, $h_j^m$, can be written as
\begin{equation}
    h_j^m = \text{ReLU}\left(\sum_{i=1}^{d_{m-1}}w_{ij}^{m-1}  h_i^{m-1} + b_j^{m-1}\right), \label{eq:activation}
\end{equation}
where $w_{ij}^m$ is the element at the $j$-th row and $i$-th column of $W^{m-1}$ and $b_j^{m-1}$ is the $j$-th index of $\vect b^{m-1}$.  
To model ReLU in a MIP for the $j$-th unit in the $m$-th layer, we use the variables $\hat h_j^{m}$, $\check h_j^m$ and $\hat h_i^{m-1}$ for $i = 1,\ldots, d_{m-1}$.  
The ReLU activation is then modeled with the following constraints:
\begin{subequations}
\begin{align}
    \sum_{i=1}^{d_{m-1} }w_{ij}^{m-1} \hat h_i^{m-1} + b_j^{m-1} & = \hat h_j^m - \check h_j^m, \label{eq:re_hu} \\
    z_j^m = 1 \Rightarrow \hat h_j^m \le 0, & \label{eq:re_hz}\\
    z_j^m = 0 \Rightarrow \check h_j^m \le 0, & \label{eq:re_sz}\\
    \hat h_j^m, \check h_j^m \ge 0, & \label{eq:re_hs}\\
    z_j^m \in \{0,1\}, & \label{eq:re_z}
\end{align}
\end{subequations}
where the logical constraints in \Cref{eq:re_hz} and \Cref{eq:re_sz} are translated into big-M constraints by MIP solvers.  
To verify the correctness of this formulation, observe that constraints \eqref{eq:re_hz} and \eqref{eq:re_sz} in conjunction with the fact the binary $z_j^m$ ensures that at most one of $\hat h_j^m$ and $\check h_j^m$ are non-zero.  Furthermore, since both $\hat h_j^m$ and $\check h_j^m$ are non-negative, if $\sum_{i=1}^{d_{m-1} }w_{ij}^{m-1} \hat h_i^{m-1} + b_j^{m-1} > 0$, then it follows that $\hat h_j^m>0$ and $\check h_j^m=0$.  
If negative, then $\hat h_j^m=0$ and $\check h_j^m>0$.  
Thus, we have that if the left-hand side of \eqref{eq:re_hu} is positive, $\hat h_j^m$ will be positive; if it is negative, then $\hat h_j^m=0$; this is an exact representation of the ReLU function.

\section{Related Work}
\label{sec:related_work}

\subsection{Machine Learning for Nested Optimization}
Machine learning has recently been employed to solve nested optimization problems; by ``nested", we mean optimization problems whose objective or constraints involve another optimization. For example, \cite{nair2018learning, shen2021learning, jiang21learn2defend, xiong2020adversarial, shao2022_learning_ro} directly predict a binary or continuous solution vector.  The major limitation with predicting solutions directly is the inability to handle variable integrality and hard constraints.  In addition, only \cite{nair2018learning} consider 2SP, whereas the others focus on bi-level problems with a single inner optimization, rather than the expectation as in stochastic programming.  
For stochastic programming, there has been a significant interest in the integration of learning to enhance prevalent solution techniques. We specifically discuss three areas of related work: learning-enabled optimization, learning-based algorithms for stochastic programming, and scenario reduction for stochastic programming. 

The line of work on learning-enabled optimization \citep{deng2022predictive,liu2022coupled,diao2020distribution} introduced ``predictive stochastic programming" to leverage contextual information when formulating SP models.
This is in contrast to our approach, which leverages predictions to reduce computing times in a non-contextual 2SP setting.  That being said,~\method{} admits extensions to the contextual setting by including the context information during training.

In recent years, several studies have explored the use of integrating predictions within stochastic programming algorithms for computational improvements.  
\citet{donti2017task} proposed an end-to-end approach to directly optimize a task-loss for contextual stochastic programming problems by differentiating through the argmin operator, specifically for strongly convex problems.
\citet{dai2022neural} developed a model to solve  multi-stage linear optimization problems by learning the piece-wise value function of the nested problems.  
\citet{larsen2022fast} leveraged predictions to improve an exact decomposition-based algorithm for 2SP. 
\method{} differs from these approaches as it can be applied to problems with both hard constraints and integer/non-linear second-stage problems.

Lastly, another related research direction for learning-based stochastic programming is scenario reduction, which reduces the complexity of the stochastic programming problem by finding a smaller set of ``representative scenarios".
Many of these approaches \citep{dupavcova2003scenario, romisch2009scenario, beraldi2014clustering, prochazka2020scenario, keutchayan2021problem} perform some form of clustering to reduce the number of scenarios and then solve a smaller surrogate problem with these scenarios.
Recently, \citet{wu2022learning} used a conditional variational autoencoder to learn scenario embeddings and perform clustering on them for scenario reduction.  To find representative scenarios, they use semi-supervised learning with respect to the second-stage cost. However, these predictions are not leveraged explicitly in the optimization as is done with~\method{}.
\citet{bengio2020learning} predicts a representative scenario for an input scenario set and use it to form a smaller surrogate problem.
They show that using the predicted representative scenario, a near-optimal first-stage decision can be obtained by solving the surrogate. However, their method requires some domain expertise as it relies on the problem structure to build the representative scenario for training.

\subsection{Neural Network Embeddings}


\method{} can be broadly classified as both a learning-based scenario reduction approach and a learning-accelerated heuristic for stochastic programming.  The reason for this is that \method{} reduces the computational complexity introduced by the scenarios by computing a compact representation that is then leveraged to formulate an approximation to the EF.  
We specifically leverage the recent line of work by \citet{cheng2017maximum, tjeng2017evaluating, fischetti2018deep, serra2018bounding}, which studies the problem of embedding a trained neural network with ReLU activation into a MIP.
The works of \citet{anderson2020strong} and \citet{grimstad2019relu} present MIP encoding formulations with tighter LP relaxations by appropriately setting the big-$M$ constraints, leading to reduced solving time.
The growing interest in embedding predictive models in MIPs has led to the development of libraries such as \texttt{JANOS} \citep{bergman2022janos} and OMLT \citep{ceccon2022omlt}. 
\cite{say2017nonlinear,grimstad2019relu,murzakhanov2020neural,KATZ202011350,KODY2022108282} propose the use of embedded neural networks to formulate surrogate MIPs for intractable and non-linear constraints in optimization problems.  To the best of our knowledge, \method{} is the first approach that employs this technique in stochastic programming or more generally for the simplification of nested optimization problems. 

\section{The~\method{} Framework}
\label{sec:methodology}



In this section, we present two neural architectures, the corresponding surrogate problems that approximate a given~\twostocprog{}, and a data collection strategy. 
Figure~\ref{fig:block-diagram} summarizes the~\method{} framework. 



\subsection{Neural Network Architectures}
\label{sec:methodology_lm}

\begin{figure}
    \centering
    \includegraphics[width=\linewidth]{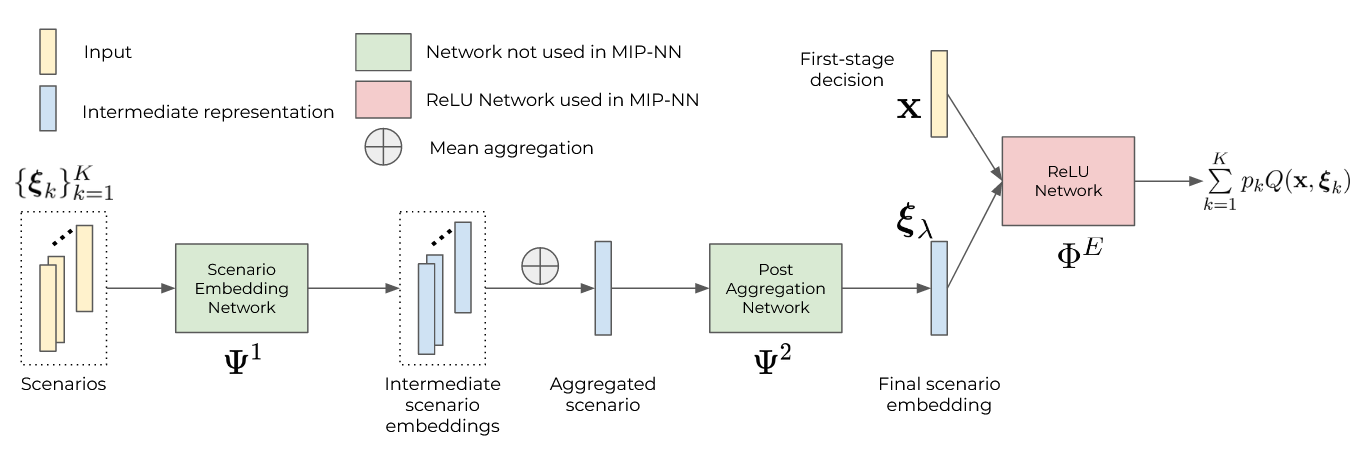}
    \caption{\EVA{} architecture diagram.}
    \label{fig:nsp_single_cut_2}
\end{figure}

We propose two distinct neural architectures for predicting the second-stage costs:
\EVA{} approximates the expected value of the second-stage cost of \textit{a set of scenarios}, whereas~\PSA{} approximates the \textit{per-scenario} value of the second-stage cost for \textit{a single scenario}.
\vspace{-.2cm}
\paragraph{\EVA}
(\Cref{fig:nsp_single_cut_2}) learns a mapping from $\big(\vect x, \{\vects \xi_k\}_{k=1}^K\big) \rightarrow \sum_{k=1}^{K} p_k Q(\vect x, \vects \xi_k)$. In words, the model takes in a first-stage solution $\vect x$ and any finite set of scenarios sampled from $\Xi$, and outputs a prediction of the expected second-stage objective value.
We embed the scenario set $\{\vects \xi_k\}_{k=1}^K$ into a latent space by passing each scenario, independently, through the same neural network $\Psi^1$, then performing mean-aggregation over the resulting $K$ embeddings.
The aggregated embedding is passed through another network, $\Psi^2$, to obtain the final embedding of the scenario set, $\xi_\lambda$.
This embedding, representing the scenario set à-la-DeepSets~\citep{deepsets17}, is appended to the input first-stage decision and passed through a ReLU feed-forward network $\Phi^{E}$ to predict the expected second-stage value.
Hence, the final output is such that $\Phi^{E}(\vect x, \Psi^2(\oplus_{k=1}^K \Psi^1(p_k, \vects \xi_k))) \approx \sum_{k=1}^{K} p_k Q(\vect x, \vects \xi_k)$.
Note that the embedding networks, $\Psi^1$ and $\Psi^2$, can be arbitrarily complex as only the latent representation is embedded into the approximate MIP. Also, although $\Psi^1$ is trained using $K$ scenarios, once the networks are trained, they can be used with any (potentially much larger) finite number of scenarios.
\vspace{-.2cm}
\paragraph{\PSA{}}
learns a mapping $\Phi^{P}$ from $(\vect x, \vects \xi) \rightarrow Q(\vect x, \vects \xi)$ for $\vects \xi$ sampled from $\Xi$.
Once the mapping $\Phi^{P}$ is learned, we can approximate the expected second-stage objective value for any finite set of scenarios as $\sum_{k=1}^{K} p_k Q(\vect x, \vects \xi_k) \approx \sum_{k=1}^{K} p_k \Phi^{P}(\vect x, \vects \xi_k)$.  
$\Phi^{P}$ is a feed-forward neural network with input given by the concatenation of $\vect x$ and $\vects \xi$.

\subsection{Neural Network Embedding for \twostocprog{}}
\label{sec:methodology_embed}

We now describe the surrogate MIP for both the \EVA{} and \PSA{} learning models from the preceding section.
Let $\Lambda$ represent the number of predictions made by the neural network.
For the \EVA{} case, $\Lambda=1$ as we only predict the expected second-stage value for a set of scenarios. 
In the \PSA{} case, $\Lambda=K$ as we predict the second-stage value for each scenario. In this section, we use $[M]$ to denote $\{1, \dots, M\}$ for $M\in\mathbb{Z}_+$. 

Let $\hat h_{j}^{m,\lambda}$ represent the ReLU output for the $j$-th hidden unit in the $m$-th hidden layer for output $\lambda$, for all $m\in[\ell - 1]$, $j\in [d_m]$, and $\lambda\in[\Lambda]$. 
Suppose $\check h_{j}^{m,\lambda}$ is a slack variable used to model the $ReLU$ output for the $j$-th hidden unit in the $m$-th hidden layer for scenario $k$, for all $m\in[\ell - 1]$, $j\in [d_m]$, and $\lambda\in[\Lambda]$.
Let $z_{j}^{m,\lambda}$ be a binary variable used to ensure that at most one of $\hat h_{j}^{m,k}$ and $\check h_{j}^{m,k}$ are non-zero.  This variable is defined for all $m\in[\ell - 1]$, $j\in [d_m]$, and $\lambda\in[\Lambda]$.
Suppose $\NNoutput_\lambda$ is the $\lambda$-th prediction by the neural network, for all $\lambda\in[\Lambda]$.

With the above variables we can define an approximation to EF as given in \cref{eq:relu_mip}. 
The objective function \eqref{eq:relu_obj} minimizes the sum of the cost of the first-stage decisions and the approximate cost of the second-stage value.  Constraints \eqref{eq:relu_in}-\eqref{eq:relu_out} propagate a first-stage solution $\vect x$ to the output of the neural network for each scenario. Constraints \eqref{eq:relu_h_z}-\eqref{eq:relu_h_s} ensure the prediction of the neural network is respected.  Constraint \eqref{eq:relu_x} ensures the feasibility of the first-stage solution.

In this approximation, we introduce a number of additional variables and big-M constraints.  Specifically, for a neural network with $H$ hidden units, we introduce $\Lambda \cdot H$ additional binary variables for $z_{j}^{m,\lambda}$.  In addition, we introduce $2\cdot \Lambda \cdot H$ continuous variables for $\hat h_{j}^{m,\lambda}$ and $\check h_{j}^{m,\lambda}$.  Lastly, we require  an additional $\Lambda$ variables for the output of the network.  
Although the number of variables we introduce in this approximation is quite large, we hypothesize that the resulting MIP will be easier to solve than the extensive form, in particular, when the second-stage problem is nonlinear. 
In the remainder of the paper, we refer to the surrogate MIP given in \eqref{eq:relu_mip} as MIP-NN.
\begin{subequations}
\begin{align}
    \min & \quad \vect c^\intercal \vect x + \sum_{\lambda=1}^{\Lambda} p_\lambda \NNoutput_\lambda \label{eq:relu_obj}\\
    \text{s.t.} & \quad \sum_{i=1}^{d_0} w_{ij}^{0} [\vect x, \vects \xi_\lambda]_i + b_j^{0} = \hat h_{j}^{1,\lambda} - \check h_{j}^{1,\lambda} & \quad \forall \hspace{3pt} j\in[d_1], \lambda\in[\Lambda], 
    \label{eq:relu_in} \\
    & \quad \sum_{i=1}^{d_{m-1}} w_{ij}^{m-1} \hat h_{i}^{m-1, \lambda} + b_j^{m-1} = \hat h_j^{m,\lambda} - \check h_{j}^{m,\lambda} & \quad \forall \hspace{3pt} m\in[\ell - 1], j\in [d_m], \lambda\in[\Lambda], \label{eq:relu_hl}\\
    & \quad \sum_{i=1}^{d_{\ell}} w_{ij}^{\ell} \hat h_i^{\ell,\lambda} + b^{\ell} \le  \NNoutput_\lambda & \quad \forall \lambda\in[\Lambda], 
    \label{eq:relu_out}\\
    & \quad z_{j}^{m,\lambda} = 1 \Rightarrow \hat h_{j}^{m,\lambda} = 0 & \quad  \forall \hspace{3pt} m\in[\ell - 1], j\in [d_m], \lambda\in[\Lambda], \label{eq:relu_h_z}\\
    & \quad {z_{j}^{m,\lambda}} = 0 \Rightarrow \check h_{j}^{m,\lambda} = 0 & \quad  \forall \hspace{3pt} m\in[\ell - 1], j\in [d_m], \lambda\in[\Lambda] \label{eq:relu_s_z},\\
    & \quad {z_{j}^{m,\lambda} \in \{0,1\}} & {\quad  \forall \hspace{3pt} m\in[\ell - 1], j\in [d_m], \lambda\in[\Lambda]}
    {\label{eq:relu_z}},\\
    & \quad {\hat h_{j}^{m,\lambda}, \check h_{j}^{m,\lambda} \ge 0} & {\quad  \forall \hspace{3pt} m\in[\ell - 1], j\in [d_m], \lambda\in[\Lambda]}
    {\label{eq:relu_h_s}},\\
    & \quad { \vect x \in \mathcal{X}} 
    {\label{eq:relu_x}}
\end{align}
\label{eq:relu_mip}
\end{subequations}

\subsection{Data Generation}

\label{sec:methodology_dg}

A diverse dataset of input-output pairs is needed to train~\method{}'s supervised second-stage value approximation.
To generate such a dataset for a given \twostocprog{} problem, we adopt an iterative procedure.
We begin by generating a random feasible first-stage decision.  For the \EVA{} case, we sample a set of scenarios with random cardinality $K^{\prime} $ from the uncertainty distribution.  Here, $K^{\prime}$ should be chosen to balance the trade-off between the time spent to generate a sample of second-stage values for a given first-stage solution and the time to estimate the expected second-stage value for a set of first-stage decisions in a given time budget.  Specifically, if $K^{\prime}$ is large, then on average more time will be spent in determining the expected value using a large number of scenarios, while for a small $K^{\prime}$, the first-stage decision space will be explored more since expected value estimates would be obtained faster.  
For a given input, i.e., a first-stage decision and set of scenarios, we then compute a label by calculating the expected second-stage value $\sum_{k^{\prime}=1}^{K^{\prime}}p_{k^{\prime}} Q_{k^{\prime}}(\cdot, \vects \xi_{k^{\prime}})$.

For the \PSA{} case, at each iteration of the data generation procedure, we sample a single scenario from the uncertainty distribution.
For a given input of a first-stage decision and scenario we generate a label by calculating its second-stage value $Q(\cdot, \cdot)$.  Last, the input-output pair is added to the dataset. 

This data generation procedure is fully parallelizable over the second-stage problems to be solved.

\subsection{\EVA{} vs. \PSA{} Trade-offs}

The \EVA{} and \PSA{} architectures exhibit trade-offs in terms of the learning task and the resulting surrogate optimization problem. 

\paragraph{Training.}
In data collection, both models require solving second-stage problems with a fixed first-stage solution to obtain the label. A sample in for \PSA{} requires solving only a single optimization problem, whereas a sample for \EVA{} requires solving at most $K^{\prime}$ second-stage problems. 
As this process is offline and highly parallelizable, this trade-off is easy to mitigate.
As for training, \EVA{} operates on a subset of scenarios which makes for an exponentially larger input space. 
Despite the large input space, our experiments show that the \EVA{} model in the training converges quite well and in many cases the embedded model outperforms the \PSA{} model.  

\paragraph{Surrogate Optimization Problem.}
As the ultimate goal is embedding the trained model into a MIP, the trade-off in this regard becomes quite important.  Specifically, for $K$ scenarios, the \PSA{} model will have $K$ times more binary and continuous variables than the \EVA{} model.
For problems with a large number of scenarios, this makes the \EVA{} model much more appealing, smaller and likely faster to solve. Furthermore, it allows for much larger networks given that only a single copy of the network is embedded.

\section{Experimental Setup}
\label{sec:setup}

All experiments were run on a computing cluster with an Intel Xeon CPU E5-2683 and Nvidia Tesla P100 GPU with 64GB of RAM (for training).  Gurobi 9.1.2 \citep{gurobi} was used as the MIP solver.  Scikit-learn 1.0.1 \citep{scikit-learn} and Pytorch 1.10.0 \citep{NEURIPS2019_9015} were used for supervised learning.  The code to reproduce all of the experiments is available at 
\url{https://github.com/khalil-research/Neur2SP}.

\textbf{2SP Problems}:
We evaluate our approach on four 2SP problems that are commonly considered in the literature: a two-stage stochastic variant of the Capacitated Facility Location Problem (CFLP)~\citep{cornuejols1991comparison}, an Investment Problem (INVP)~\citep{schultz1998solving}, the Stochastic Server Location Problem (SSLP) \citep{ntaimo2005million}, and the Pooling Problem (PP) \citep{audet2004pooling, gupte2017relaxations, gounaris2009computational, haverly1978studies}. \cref{tab:problem_summary} summarizes the types of first and second-stage variables for these problems and Appendix~\ref{app:spp} includes their detailed descriptions.

\begin{table}[t]
\centering
\resizebox{0.75\textwidth}{!}{
\begin{tabular}{llllll}
    \toprule
    {Problem} & {First stage } & {Second Stage } & {Objective} & {Constraints}  & {Objective Sense}\\
    \midrule
    CFLP & Binary       & Binary             &  Linear          & Linear        & Minimization \\
    INVP & Continuous   & Binary             &  Linear          & Linear        & Minimization \\
    SSLP & Binary       & Binary             &  Linear          & Linear        & Minimization\\
    PP   & Binary       & Continuous         &  Bilinear        & Bilinear      & Maximization \\
    \bottomrule
\end{tabular}}
\caption{Problem class characteristics.}
\label{tab:problem_summary}
\end{table}

\textbf{Baselines}:
We consider two baselines.  The first is \extform{}, which is perhaps the only generic approach that can be applied for both integer and nonlinear second-stage problems.  
We limit the solving time of \extform{} to
3 hours.
Additionally, we compare against an embedded trained linear regressor rather than a neural network, but defer these results to Appendix~\ref{app:objective} as the solution quality is quite poor in comparison to the neural network models.

\textbf{Model \& Dataset Selection }:
As is common in supervised learning, model selection and the size of the training set can have a significant impact on model performance. We present detailed experiments for model selection and dataset sizing in Appendix~\ref{app:ms_ds}.  As a brief summary, we use random search over 100 hyperparameter configurations for model selection, and observe that accuracy on a validation set is rather insensitive to hyperparameter choices.  For the size of the dataset, we observe diminishing returns when increasing the size beyond 5000 samples.

\textbf{Data Generation \& Supervised training times }:
As the data generation and training can be done offline and are both parallelizable, 
we report the total times in Table~\ref{tab:offline_times} and defer more specific timing details to Appendix~\ref{app:dg_tr}.
We note that for data generation, a single sample can be obtained in less than two seconds for all instances, and in many cases much faster.
The training times are within the range of 120 to 2100 seconds.  For the \EVA{} data generation, we choose $K^{\prime}=100$, a number of scenarios which was quick to label while exposing the model to a reasonably large number of scenarios in some cases.  The combined time for data generation and model training is typically less than 3 hours (i.e., the time given to EF) and depending on the problem, may be much less.

\begin{table*}[t]\centering\resizebox{0.7\textwidth}{!}{
\begin{tabular}{l|rr|rr|rr}
\toprule
Problem & \multicolumn{2}{c}{Data Generation Time} & \multicolumn{2}{c}{Training Time} & \multicolumn{2}{c}{Total Time } \\
\cmidrule{2-7}
{} & {\EVA{}} & {\PSA{}} & {\EVA{}\ } & {\PSA{}\ } & {\EVA{}\ \ } & {\PSA{}\ \ } \\
\midrule
CFLP\_10\_10 & 1,823.07 & 13.59 & 667.28 & 127.12 & 2,490.35 & 140.71 \\
CFLP\_25\_25 & 4,148.83 & 112.83 & 2,205.23 & 840.07 & 6,354.06 & 952.90 \\
CFLP\_50\_50 & 7,697.91 & 135.57 & 463.71 & 128.11 & 8,161.62 & 263.68 \\
\midrule
SSLP\_10\_50 & 942.10 & 22.95 & 708.86 & 116.17 & 1,650.96 & 139.13 \\
SSLP\_15\_45 & 929.27 & 16.35 & 1,377.21 & 229.42 & 2,306.48 & 245.77 \\
SSLP\_5\_25 & 860.74 & 13.18 & 734.02 & 147.44 & 1,594.75 & 160.62 \\
\midrule
INVP\_B\_E & 8,951.27 & 4.17 & 344.87 & 1,000.14 & 9,296.13 & 1,004.31 \\
INVP\_B\_H & 9,207.90 & 4.22 & 1,214.54 & 607.49 & 10,422.43 & 611.71 \\
INVP\_I\_E & 8,759.83 & 4.34 & 2,115.25 & 680.93 & 10,875.08 & 685.27 \\
INVP\_I\_H & 8,944.65 & 3.32 & 393.82 & 174.26 & 9,338.47 & 177.58 \\
\midrule
PP & 1,202.11 & 14.86 & 576.08 & 367.25 & 1,778.19 & 382.11 \\
\bottomrule
\end{tabular}}
\caption{
Computing times (in seconds) for data generation and training.  Data was generated in parallel with 43 processes. 
}
\label{tab:offline_times}
\end{table*}

\section{Results \& Discussion}
\label{sec:results}


In this section, we report the results of~\method{} across the four problem settings.   As is standard in \twostocprog{}, we evaluate a single ``base'' instance across varying scenario sets and sizes. For example, in CFLP, and for a ``base'' instance with 10 facilities and 10 customers, one can generate an instance by sampling any number of scenarios. 
An important advantage of our approach is that we can apply a single trained model to an instance with an arbitrary number of scenarios.
For example, the same trained model is used for CFLP\_10\_10\_\{100,500,1000\}.  

Tables~\ref{tab:res_CFLP} through \ref{tab:res_PP} report the gaps between approaches, solving times, and the time which \extform{} takes to achieve the same solution quality as \EVA{} and \PSA{}.  In addition, we include supplementary results with the objective values in Appendix~\ref{app:objective} and non-aggregated results for the SSLP SIPLib instances in Appendix~\ref{app:siplib}.
For SSLP and CFLP, each row represents mean statistics across 11 and 10 instances generated by sampling different scenario sets of a given size, respectively.
However, for INVP and PP, each row represents the statistics across 1 instance. 
Originally, both these problems have infinite support as the uncertainty distributions are assumed to be continuous. 
To manage the complexity, these distributions are typically transformed to have finite support by uniformly sampling equidistant points over the continuous domain.
We adopt this same procedure from the literature, leading to a static set of scenarios for a given scenario set size.

\begin{table*}[t]\centering\resizebox{0.82\textwidth}{!}{
\begin{tabular}{l|rr|rrr|rrrr}
\toprule
Problem & \multicolumn{2}{c|}{Obj. Difference (\%)} & \multicolumn{3}{c|}{Solving Time} & \multicolumn{4}{c}{EF time to } \\
\cmidrule{2-10}
{} & {EF-\EVA{}} & {EF-\PSA{}} & {\EVA{}} & {\PSA{}} & {EF} & {\EVA{}\ } & {\ } & {\PSA{}\ } & {} \\
\midrule
CFLP\_10\_10\_100 & 2.58 & 1.65 & \textbf{0.38} & 8.28 & 4,410.60 & 8.87 & (0) & 12.69 & (0) \\
CFLP\_10\_10\_500 & 2.41 & 0.94 & \textbf{0.60} & 206.30 & 10,800.17 & 415.89 & (0) & 2,034.73 & (0) \\
CFLP\_10\_10\_1000 & 0.94 & -0.67 & \textbf{0.64} & 856.77 & 10,800.87 & 580.50 & (0) & 7,551.00 & (8) \\
CFLP\_25\_25\_100 & -0.75 & -0.75 & \textbf{0.44} & 4.86 & 10,800.06 & - & (10) & - & (10) \\
CFLP\_25\_25\_500 & -3.62 & -3.62 & \textbf{0.54} & 26.41 & 10,800.14 & - & (10) & - & (10) \\
CFLP\_25\_25\_1000 & -1.32 & -1.32 & \textbf{0.58} & 54.45 & 10,800.36 & - & (10) & - & (10) \\
CFLP\_50\_50\_100 & -0.43 & -1.29 & \textbf{1.66} & 21.10 & 10,800.05 & 5,637.98 & (6) & 2,334.04 & (9) \\
CFLP\_50\_50\_500 & -9.58 & -10.71 & \textbf{1.25} & 173.63 & 10,806.15 & - & (10) & - & (10) \\
CFLP\_50\_50\_1000 & -16.62 & -17.50 & \textbf{1.44} & 572.12 & 10,805.82 & - & (10) & - & (10) \\
\bottomrule
\end{tabular}}
\caption{CFLP results: each row represents an average over ten 2SP instances with varying scenario sets. 
``Obj. Difference'' for method EF-\{\EVA{}, \PSA{}\} is the percent relative objective value of \{\EVA{}, \PSA{}\} to EF; a negative (positive) value of $-g\%$ ($g\%$) indicates that \{\EVA{}, \PSA{}\} finds a solution that is $g\%$ better (worse) than EF's for the minimization problem.
``Solving Time'' is the time in which \{\EVA{}, \PSA{}, EF\} are solved to optimality.  A time of $\sim$10,800 implies that the solving limit was reached. 
``EF time to'' is the time in which EF achieves a solution of the same quality as \{\EVA{}, \PSA{}\}.  To the right in parentheses is the number of instances for which EF failed to find a solution that is as good as \{\EVA{}, \PSA{}\}.  If EF did not find any feasible solution, then the entry is left as ``-''.  
All times are in seconds.}
\label{tab:res_CFLP}
\end{table*}

\begin{table*}[t]\centering\resizebox{0.82\textwidth}{!}{
\begin{tabular}{l|rr|rrr|rrrr}
\toprule
Problem & \multicolumn{2}{c|}{Obj. Difference (\%)} & \multicolumn{3}{c|}{Solving Time} &
\multicolumn{4}{c}{EF time to } \\
\cmidrule{2-10}
{} & {EF-\EVA{}} & {EF-\PSA{}} & {\EVA{}} & {\PSA{}} & {EF} & {\EVA{}\ } & {\ } & {\PSA{}\ } & {} \\
\midrule
SSLP\_10\_50\_50 & 0.00 & 0.00 & \textbf{0.11} & 5.83 & 10,800.48 & 228.06 & (0) & 228.06 & (0) \\
SSLP\_10\_50\_100 & -0.00 & -0.00 & \textbf{0.11} & 13.09 & 10,800.21 & 145.35 & (0) & 145.35 & (0) \\
SSLP\_10\_50\_500 & -0.00 & -0.00 & \textbf{0.14} & 129.44 & 10,802.82 & 7,359.85 & (4) & 7,359.85 & (4) \\
SSLP\_10\_50\_1000 & -55.21 & -55.21 & \textbf{0.13} & 466.38 & 10,800.47 & - & (11) & - & (11) \\
SSLP\_10\_50\_2000 & -102.69 & -102.69 & \textbf{0.14} & 2,182.31 & 10,800.17 & - & (11) & - & (11) \\
SSLP\_15\_45\_5 & 3.10 & 18.71 & \textbf{0.32} & 0.34 & 2.54 & 0.75 & (0) & 0.12 & (0) \\
SSLP\_15\_45\_10 & 2.98 & 18.47 & \textbf{0.31} & 0.58 & 1,976.62 & 2.72 & (0) & 0.20 & (0) \\
SSLP\_15\_45\_15 & 2.53 & 18.90 & \textbf{0.33} & 0.86 & 2,052.76 & 1.84 & (0) & 0.34 & (0) \\
SSLP\_5\_25\_50 & 0.12 & 1.78 & \textbf{0.20} & 1.14 & 2.24 & 1.94 & (0) & 1.97 & (0) \\
SSLP\_5\_25\_100 & 0.02 & 1.60 & \textbf{0.18} & 1.83 & 8.43 & 8.04 & (0) & 7.75 & (1) \\
\bottomrule
\end{tabular}}
\caption{SSLP results: each row represents an average over eleven 2SP instances with varying scenario sets, one of which being the instance from \cite{ahmed2015siplib}. Columns are as in Table~\ref{tab:res_CFLP}.}
\label{tab:res_SSLP}
\end{table*}

\begin{table*}[t]\centering\resizebox{0.75\textwidth}{!}{
\begin{tabular}{l|rr|rrr|rr}
\toprule
Problem & \multicolumn{2}{c|}{Obj. Difference (\%)} & \multicolumn{3}{c|}{Solving Time} & \multicolumn{2}{c}{EF time to } \\
\cmidrule{2-8}
{} & {EF-\EVA{}} & {EF-\PSA{}} & {\EVA{}} & {\PSA{}} & {EF} & {\EVA{}\ } & {\PSA{}\ } \\
\midrule
INVP\_B\_E\_4 & 9.54 & 3.01 & 0.36 & 0.34 & \textbf{0.02} & 0.02 & 0.02 \\
INVP\_B\_E\_9 & 7.54 & 2.00 & 0.31 & 0.53 & \textbf{0.04} & 0.03 & 0.03 \\
INVP\_B\_E\_36 & 2.72 & 4.96 & 0.30 & 9.53 & \textbf{0.08} & 0.02 & 0.02 \\
INVP\_B\_E\_121 & 1.37 & 2.42 & \textbf{0.33} & 86.42 & 1.69 & 0.06 & 0.02 \\
INVP\_B\_E\_441 & 2.80 & 2.43 & \textbf{0.37} & 4,342.19 & 117.59 & 0.78 & 1.15 \\
INVP\_B\_E\_1681 & 1.36 & - & \textbf{0.34} & - & 10,800.01 & 17.41 & 0.00 \\
INVP\_B\_E\_10000 & -1.48 & - & \textbf{0.36} & - & 10,803.98 & - & 0.00 \\
\hline
INVP\_B\_H\_4 & 8.81 & 9.50 & 0.46 & 0.25 & \textbf{0.01} & 0.01 & 0.01 \\
INVP\_B\_H\_9 & 5.04 & 5.04 & 0.30 & 0.57 & \textbf{0.03} & 0.02 & 0.02 \\
INVP\_B\_H\_36 & 1.61 & 1.61 & \textbf{0.26} & 6.79 & 1.29 & 0.01 & 0.01 \\
INVP\_B\_H\_121 & 1.77 & 1.77 & \textbf{0.33} & 45.89 & 34.69 & 0.01 & 0.01 \\
INVP\_B\_H\_441 & 2.13 & 5.50 & \textbf{0.28} & 1,870.42 & 217.46 & 2.21 & 0.21 \\
INVP\_B\_H\_1681 & -0.71 & - & \textbf{0.36} & - & 10,800.01 & - & 0.00 \\
INVP\_B\_H\_10000 & -2.72 & - & \textbf{0.36} & - & 10,800.03 & - & 0.00 \\
\hline
INVP\_I\_E\_4 & 12.83 & 0.00 & 0.38 & 0.23 & \textbf{0.01} & 0.01 & 0.01 \\
INVP\_I\_E\_9 & 7.40 & 2.64 & 0.27 & 0.35 & \textbf{0.06} & 0.01 & 0.02 \\
INVP\_I\_E\_36 & 5.48 & 5.17 & 0.27 & 1.39 & \textbf{0.04} & 0.01 & 0.01 \\
INVP\_I\_E\_121 & 5.30 & 4.49 & \textbf{0.29} & 49.51 & 1.65 & 0.02 & 0.03 \\
INVP\_I\_E\_441 & 3.00 & 0.68 & \textbf{0.26} & 2,049.93 & 46.92 & 0.08 & 0.10 \\
INVP\_I\_E\_1681 & 1.31 & 3.08 & \textbf{0.26} & 10,834.53 & 10,800.00 & 0.41 & 0.41 \\
INVP\_I\_E\_10000 & -1.35 & - & \textbf{0.30} & - & 10,800.10 & - & 0.00 \\
\hline
INVP\_I\_H\_4 & 13.78 & 12.16 & 0.35 & 0.21 & \textbf{0.02} & 0.01 & 0.01 \\
INVP\_I\_H\_9 & 9.12 & 0.81 & 0.37 & 0.31 & \textbf{0.03} & 0.01 & 0.02 \\
INVP\_I\_H\_36 & 4.97 & 3.44 & \textbf{0.36} & 1.99 & 1.27 & 0.03 & 0.03 \\
INVP\_I\_H\_121 & 4.01 & 4.99 & \textbf{0.32} & 23.10 & 7.43 & 0.07 & 0.07 \\
INVP\_I\_H\_441 & 3.15 & 3.15 & \textbf{0.32} & 1,231.48 & 10,800.00 & 0.33 & 0.33 \\
INVP\_I\_H\_1681 & -0.34 & 0.11 & \textbf{0.33} & 10,816.89 & 10,800.03 & - & 252.70 \\
INVP\_I\_H\_10000 & -1.60 & - & \textbf{0.38} & - & 10,802.10 & - & 0.00 \\
\bottomrule
\end{tabular}}
\caption{INVP results: each row represents a single instances. Columns are as in Table~\ref{tab:res_CFLP}.}
\label{tab:res_INVP}
\end{table*}

\begin{table*}[t]\centering\resizebox{0.75\textwidth}{!}{
\begin{tabular}{l|rr|rrr|rr}
\toprule
Problem & \multicolumn{2}{c|}{Obj. Difference (\%)
} & \multicolumn{3}{c|}{Solving Time} & \multicolumn{2}{c}{EF time to } \\
\cmidrule{2-8}
{} & {EF-\EVA{}} & {EF-\PSA{}} & {\EVA{}} & {\PSA{}} & {EF} & {\EVA{}\ } & {\PSA{}\ } \\
\midrule
PP\_125 & 3.25 & 36.64 & \textbf{1.51} & 144.08 & 10,800.00 & 10,717.05 & 2.48 \\
PP\_216 & 9.06 & 40.14 & \textbf{1.47} & 254.94 & 364.98 & 59.79 & 3.59 \\
PP\_343 & 0.67 & 40.85 & \textbf{1.46} & 570.64 & 10,800.00 & 1,450.54 & 9.12 \\
PP\_512 & 8.69 & 39.77 & \textbf{1.60} & 1,200.37 & 10,800.01 & 167.77 & 13.80 \\
PP\_729 & 1.38 & 37.98 & \textbf{1.62} & 3,440.19 & 10,800.01 & 5,867.67 & 36.34 \\
PP\_1000 & 5.92 & 41.32 & \textbf{1.49} & 10,853.59 & 10,800.00 & 1,596.22 & 210.48 \\
\bottomrule
\end{tabular}}
\caption{PP results: each row represents a single instances. Columns are as in Table~\ref{tab:res_CFLP}.}
\label{tab:res_PP}
\end{table*}

\subsection{Discussion}
Tables~\ref{tab:res_CFLP}--\ref{tab:res_PP} show that \EVA{} is significantly faster than other approaches, with a minimum and maximum solving time of 0.11s and 1.66s respectively, across all problems.
This highlights the scalability of the \EVA{} in terms of problem size and type, which is expected as the size of the resulting MIP is independent of the number of scenarios. 
Also, the objective difference is less than 5\% in most cases, with a minimum of -102\% and a maximum of 13.78\%.
These differences are inversely proportional to the scenario set size, which indicates that the \EVA{} is able to generalize on larger scenario sets, even though it was trained with a maximum of 100 scenarios per data point.
EF takes significantly longer to reach a solution quality similar to \EVA{}, often on the order of minutes to 3 hours. For many larger CFLP and SSLP instances, EF finds worse solutions than \EVA{} even after 3 hours. Not only is \EVA{} as good or better in solution quality, but also orders of magnitude faster. 

For the \PSA{}, we can observe that the solving time is directly proportional to the size of the problem. 
However, for the largest INVP instances, it times out without even generating a feasible solution (Table~\ref{tab:res_INVP}).
This is expected as we need to embed the trained neural network once per scenario, limiting scalability.
In terms of objective differences, we can observe that the difference improves with the increase in instance size for CFLP and SSLP, whereas
 no clear trend is visible for For INVP. However, the objective differences do not exceed 5\% in most cases. 
For PP, the objective difference is around 40\%, indicating that the \PSA{} is not able to generalize, whereas \EVA{} performs very well.  One important advantage for \PSA{} over \EVA{} is the fact that the time required for data generation and training is notably less.  In settings with limited parallel computing resources or time \PSA{} may be a more appropriate choice.


\section{Conclusion}
\label{sec:conclusion}

Two-stage stochastic programming is a powerful modeling framework for decision-making  under uncertainty. These problems are hard to solve in practice, especially when the second-stage problem is a MIP or NLP.  Finding good feasible solutions quickly thus becomes extremely important.

To that end, we proposed \texttt{Neur2SP}, a learning-based, general, and structure-agnostic approach which approximates the second-stage value function to form an easy-to-solve surrogate problem. The four problem classes we have tackled are (1) widely used in the literature, (2) vary in the types of first and second-stage problems, and (3) span a wide range in terms of number of variables, constraints, and scenarios. Through our experiments, we show that \texttt{Neur2SP} achieves high-quality solutions quickly, especially for larger instances. In 1--2 seconds, a model trained in the \texttt{Neur2SP} framework can find solutions of the same or better quality than the most generic method in the literature, EF, with the latter taking minutes to hours.

In terms of future work, this methodology can be extended in many directions.  Further innovations in the \EVA{} model architecture may improve our already positive results. 
Another direction is the extension of the general idea of embedding trained models into other complex optimization problems, such as bilevel optimization or multi-stage stochastic programming.  

Another direction for future work is a more comprehensive comparison of \method{} with algorithms that are specialized to a given problem class. However, we note that on SSLP instances, the computing times of progressive hedging \citep{rockafellar1991scenarios}, a widely used heuristic for 2SP, is on the order of hours \citep{torres2022review}. These experiments are not directly comparable as they were run on different hardware. However, this would not meaningfully impact the several order of magnitude reduction in solving time achieved by our approach. 


Lastly, in this work, we propose \EVA{} and \PSA{}, however, a natural middle ground between these models is a clustering approach which embeds a trained model for a subset of scenarios, rather than a single or the entire scenario set at evaluation time.

\noindent{\textbf{Acknowledgments:}}
Bodur would like to acknowledge support from an NSERC Discovery Grant. Dumouchelle, Patel, and Khalil acknowledge support from the Scale AI Research Chair Program and an NSERC Discovery Grant.

\clearpage
\bibliographystyle{unsrtnat}
\bibliography{references}
\clearpage

\section*{Checklist}

The checklist follows the references.  Please
read the checklist guidelines carefully for information on how to answer these
questions.  For each question, change the default \answerTODO{} to \answerYes{},
\answerNo{}, or \answerNA{}.  You are strongly encouraged to include a {\bf
justification to your answer}, either by referencing the appropriate section of
your paper or providing a brief inline description.  For example:
\begin{itemize}
  \item Did you include the license to the code and datasets? \answerYes{See Section~\ref{gen_inst}.}
  \item Did you include the license to the code and datasets? \answerNo{The code and the data are proprietary.}
  \item Did you include the license to the code and datasets? \answerNA{}
\end{itemize}
Please do not modify the questions and only use the provided macros for your
answers.  Note that the Checklist section does not count towards the page
limit.  In your paper, please delete this instructions block and only keep the
Checklist section heading above along with the questions/answers below.

\begin{enumerate}

\item For all authors...
\begin{enumerate}
  \item Do the main claims made in the abstract and introduction accurately reflect the paper's contributions and scope?
    \answerYes{}
  \item Did you describe the limitations of your work?
    \answerYes{}
  \item Did you discuss any potential negative societal impacts of your work?
    \answerNA{}
  \item Have you read the ethics review guidelines and ensured that your paper conforms to them?
    \answerYes{}
\end{enumerate}

\item If you are including theoretical results...
\begin{enumerate}
  \item Did you state the full set of assumptions of all theoretical results?
    \answerNA{}
        \item Did you include complete proofs of all theoretical results?
    \answerNA{}
\end{enumerate}

\item If you ran experiments...
\begin{enumerate}
  \item Did you include the code, data, and instructions needed to reproduce the main experimental results (either in the supplemental material or as a URL)?
    \answerYes{ All code to reproduce the experiments is available at \url{https://anonymous.4open.science/r/neural_stochastic_programming-437E}.  }
  \item Did you specify all the training details (e.g., data splits, hyperparameters, how they were chosen)?
    \answerYes{}
        \item Did you report error bars (e.g., with respect to the random seed after running experiments multiple times)?
    \answerYes{The both training and evaluation are done over varying configurations and realizations of randomness.}
        \item Did you include the total amount of compute and the type of resources used (e.g., type of GPUs, internal cluster, or cloud provider)?
    \answerYes{These are detailed in Section~\ref{sec:setup} and the Appendix.}
\end{enumerate}

\item If you are using existing assets (e.g., code, data, models) or curating/releasing new assets...
\begin{enumerate}
  \item If your work uses existing assets, did you cite the creators?
    \answerYes{All software used in the development has been cited}
  \item Did you mention the license of the assets?
    \answerNA{}
  \item Did you include any new assets either in the supplemental material or as a URL?
    \answerNA{}
  \item Did you discuss whether and how consent was obtained from people whose data you're using/curating?
    \answerNA{}
  \item Did you discuss whether the data you are using/curating contains personally identifiable information or offensive content?
    \answerNA{}
\end{enumerate}

\item If you used crowdsourcing or conducted research with human subjects...
\begin{enumerate}
  \item Did you include the full text of instructions given to participants and screenshots, if applicable?
    \answerNA{}
  \item Did you describe any potential participant risks, with links to Institutional Review Board (IRB) approvals, if applicable?
    \answerNA{}
  \item Did you include the estimated hourly wage paid to participants and the total amount spent on participant compensation?
    \answerNA{}
\end{enumerate}

\end{enumerate}


\appendix
\section{Symbols}
\label{app:symbols}
List of symbols used in the paper with their brief description.
\begin{table}[h]
    \centering
    \small
    \begin{tabular}{c|l}
        \toprule
        \multicolumn{2}{c}{Two-stage stochastic program}\\
        \midrule
        $\vect x$ & First-stage decision vector \\
        $\vect c$ & First-stage objective coefficient vector\\
        $K$ & \extform{} scenario set size\\
        $k$ & Scenario index\\
        $\vects \xi_k$ & $k^{th}$ scenario realization \\ 
        $p_k$ & Probability of scenario $k$ \\ 
        $n$ & Dimension of $\vect x$\\
        $Q(\vect x, \vects \xi)$ & Second-stage sub-problem for first-stage decision $\vect x$ and scenario $\vects \xi$\\
        $F(\vect x, \vects \xi)$ & Second-stage cost for first-stage decision $\vect x$ and scenario $\vects \xi$ \\
        $\mathcal X$ & Constraint set exclusively on the first-stage decision\\
        $\mathcal Y(\vect x, \vects \xi)$ & Scenario-specific constraint set for first-stage decision $\vect x$ and scenario $\vects \xi$\\
        \midrule
        \multicolumn{2}{c}{Neural network}\\
        \midrule
        $\ell$ & Number of layers in the network \\
        $m$ & Index over the neural network layers\\
        $d^{0}$ & Dimensionality of input layer\\
        $d^{m}$ & Dimensionality of layer $m$\\
        \NNinput & Input to the neural network\\
        \NNoutput & Output of the neural network\\
        $W$ & Weight matrix \\
        $\vect b$ & Bias \\
        $\sigma$ & Activation function \\
        $\vect h^{m}$ & $m^{th}$ hidden layer\\
        $i$ & Index over the column of weight matrix\\
        $j$ & Index over the row of weight matrix\\
        $\Phi^1$ & Scenario-encoding network\\
        $\Phi^2$ & Post scenario-aggregation network\\
        $\Psi^{E}$ & Scenario-embedding network for \EVA{} \\
        $\Psi^{P}$ & Scenario-embedding network for \PSA{} \\
        \midrule
        \multicolumn{2}{c}{MIP-NN}\\
        \midrule
        $\hat h$ & Non-negative $ReLU$ input\\
        $\check h$ & Negative $ReLU$ input\\
        $z$& Indicator variables\\
        $\Lambda$ & Number of predictions used in embedding\\
        $[M]$   &  The set $\{1,\ldots,M\}$ for an $M\in\mathbb{Z}_+$\\
        \bottomrule
    \end{tabular}
    \caption{Symbols summary}
    \label{tab:symbols}
\end{table}

\section{Stochastic Programming Problems}
\label{app:spp}
\subsection{Capacitated Facility Location (CFLP)}
\label{app:spp_cflp}
The CFLP is a decision-making problem in which a set of facility opening decisions must be made in order to meet the demand of a set of customers.  Typically this is formulated as a minimization problem, where the amount of customer demand satisfied by each facility cannot  exceed  its  capacity.   The  two-stage stochastic CFLP arises when facility opening decisions must be made prior to knowing the actual demand. For this problem, we generate instances following the procedure described in \citep{cornuejols1991comparison} and create a stochastic variant by simply generating the first-stage costs and capacities, then generate scenarios by sampling $K$ demand vectors using the distributions defined in \citep{cornuejols1991comparison}). To ensure relatively complete recourse, we introduce additional variables with prohibitively expensive objective costs in the case where customers cannot be served.  In the experiments a CFLP with $n$ facilities, $m$ customers, and $s$ scenarios is denoted by CFLP\_$n$\_$m$\_$s$.

\subsection{Investment Problem (INVP)}
\label{app:spp_invp}
The INVP is a 2SP problem studied in \citep{schultz1998solving}.  This 2SP has a set of continuous first-stage decisions which yield an immediate revenue. In the second stage, after a set of random variables are realized,  a set of binary decisions can be made to receive further profit.  In this work, we specifically consider the instance described in the example 7.3. of \citep{schultz1998solving}.  This problem has 2 continuous variables in the first stage with the domain $[0, 5]$, and 4 binary variables in the second stage.  The scenarios are given by two random discrete variables which are defined with equal probability over the range $[5,15]$.  Specifically, for $K$ scenarios, each random variable can take an equally spaced value in the range.  Although the number of variables is quite small, the presence of continuous first-stage decision has made this problem relevant within the context of other recent work such as the heuristic approach proposed in \citep{van2021loose}.  As a note, we reformulate the INVP as an equivalent minimization problem in the remainder of this work.  In the experiments an INVP instance is denoted by INVP\_$v$\_$t$, where $v$ indicates the type of second-stage variable 
(B for binary and I for integer) and $t$ indicates the type of technology matrix ($E$ for identity and $H$ for $[[2 / 3, 1 / 3], [1 / 3, 2 / 3]]$). 

\subsection{Stochastic Server Location Problem (SSLP)}
\label{app:spp_sslp}
The SSLP is a 2SP, where in the first stage a set of decisions are made to decide which servers should be utilized and a set of second-stage decisions assigning clients to servers. In this case, the random variables take binary values, which represent a client with a request occurring in the scenario or not. A more detailed description of the problem can be found in \citep{ntaimo2005million}.  In this work, we directly use the instances provided in SIPLIB \citep{ahmed2015siplib}.  In the experiments a SSLP with $n$ servers, $m$ clients, and $s$ scenarios is denoted by SSLP\_$n$\_$m$\_$s$.

\subsection{Pooling Problem (PP)}
\label{app:spp_pp}
The pooling problem is a well-studied problem in the field of mixed-integer nonlinear programming \cite{audet2004pooling, gupte2017relaxations, gounaris2009computational, haverly1978studies}.
    It can be formulated as a mixed-integer quadratically constrained quadratic program, making it the hardest problem class in our experiments. 
    
    We are given a directed graph, consisting of three disjoint sets of nodes, called the source, pool and terminal nodes.
    We need to produce and send some products from the source to the terminal nodes, using the given arcs, such that the product demand and quality constraints on the terminal nodes, along with the arc capacity constraints, are satisfied. 
    The pool nodes can be used to mix products with different qualities in appropriate quantities to generate a desired quality product. 
    The goal is to decide the amount of product to send on each arc such that the total profit from the operations is maximized. 
    We consider a stochastic version of the problem as described in the case study of \cite{li2011stochastic}.
    Here, in the first stage, we need to design the network by selecting nodes and arcs from the input graph, without knowing the quality of the product produced on source nodes and the exact demand on the terminal nodes.
    Once the uncertainty is revealed, in the second stage, we make the recourse decisions about the amount of product to be sent on each arc, such that demand and quality constraints on the terminal nodes are satisfied. 
    In our case, we have 16 binary variables in the first stage and 11 continuous variables per scenario in the second stage. 
    An instance of this problem is referred to as PP\_$s$, where $s$ is the number of scenarios.

\section{Data Generation \& Supervised Learning Times}
\label{app:dg_tr}

In this section, we report details of the data generation and training times for all problem settings in Tables~\ref{tab:dg_times} and \ref{tab:tr_times}, respectively.  For training, we split the \# samples into an 80\%-20\% train validation set, and select the best model on the validation set in the given number of epochs.

\begin{table*}[t]\centering\resizebox{0.9\textwidth}{!}{
\begin{tabular}{l|rrr|rrr}
\toprule
Problem & \multicolumn{3}{c}{\EVA{}} & \multicolumn{3}{c}{\PSA{}}  \\
\cmidrule{2-7}
{} & {\# samples} & {Time per sample} & {Total time} & {\# samples} & {Time per sample} & {Total time} \\
\midrule
CFLP\_10\_10 & 5,000 & 0.36 & 1,823.07 & 10,000 & 0.00 & 13.59 \\
CFLP\_25\_25 & 5,000 & 0.83 & 4,148.83 & 10,000 & 0.01 & 112.83 \\
CFLP\_50\_50 & 5,000 & 1.54 & 7,697.91 & 10,000 & 0.01 & 135.57 \\
\midrule
SSLP\_10\_50 & 5,000 & 0.19 & 942.10 & 10,000 & 0.00 & 22.95 \\
SSLP\_15\_45 & 5,000 & 0.19 & 929.27 & 10,000 & 0.00 & 16.35 \\
SSLP\_5\_25 & 5,000 & 0.17 & 860.74 & 10,000 & 0.00 & 13.18 \\
\midrule
INVP\_B\_E & 5,000 & 1.79 & 8,951.27 & 10,000 & 0.00 & 4.17 \\
INVP\_B\_H & 5,000 & 1.84 & 9,207.90 & 10,000 & 0.00 & 4.22 \\
INVP\_I\_E & 5,000 & 1.75 & 8,759.83 & 10,000 & 0.00 & 4.34 \\
INVP\_I\_H & 5,000 & 1.79 & 8,944.65 & 10,000 & 0.00 & 3.32 \\
\midrule
PP & 5,000 & 0.24 & 1,202.11 & 10,000 & 0.00 & 14.86 \\
\bottomrule
\end{tabular}}
\caption{Data generation samples and times.  Data was generated in parallel with 43 processes.  All times in seconds.}
\label{tab:dg_times}
\end{table*}


\begin{table*}[t]\centering\resizebox{0.4\textwidth}{!}{
\begin{tabular}{l|rrr}
\toprule
{} & {\EVA{}} & {\PSA{}} & {LR} \\
\midrule
CFLP\_10\_10 & 667.28 & 127.12 & 0.53 \\
CFLP\_25\_25 & 2,205.23 & 840.07 & 0.28 \\
CFLP\_50\_50 & 463.71 & 128.11 & 0.75 \\
\midrule
SSLP\_10\_50 & 708.86 & 116.17 & 0.63 \\
SSLP\_15\_45 & 1,377.21 & 229.42 & 0.57 \\
SSLP\_5\_25 & 734.02 & 147.44 & 0.05 \\
\midrule
INVP\_B\_E & 344.87 & 1,000.14 & 0.02 \\
INVP\_B\_H & 1,214.54 & 607.49 & 0.02 \\
INVP\_I\_E & 2,115.25 & 680.93 & 0.02 \\
INVP\_I\_H & 393.82 & 174.26 & 0.02 \\
\midrule
PP & 576.08 & 367.25 & 0.05 \\
\bottomrule
\end{tabular}}
\caption{Training times.  All times in seconds.}
\label{tab:tr_times}
\end{table*}

\section{Objective Results}
\label{app:objective}

In this section, we report the objective for the first-stage solutions obtained by each approximate MIP and the objective of EF (either optimal or at the end of the solving time).  In addition, we report the objective of the approximate MIP.  See Tables~\ref{tab:obj_CFLP} through \ref{tab:obj_PP} for results.  As mentioned in the main paper, the results from linear regressor (LR) are quite poor, with a significantly worse objective in almost every instance.  This is not surprising as a linear function will not likely have the capacity to estimate the integer and non-linear second-stage objectives.  For both \EVA{} and \PSA{} we can see that the true objective and the approximate-MIP objective are relatively close for all of the problem settings, further indicating that the neural network embedding is a useful approximation to the second-stage expected cost.

\begin{table*}[t]\centering\resizebox{0.9\textwidth}{!}{
\begin{tabular}{l|rrrr|rrr}
\toprule
Problem & \multicolumn{4}{c}{True objective} & \multicolumn{3}{c}{Approximate-MIP objective} \\
\cmidrule(lr){2-5}
\cmidrule(lr){6-8}
{} & {\EVA{}} & {\PSA{}} & {LR} & {EF} & {\EVA{}\ } & {\PSA{}\ } & {LR\ } \\
\midrule
CFLP\_10\_10\_100 & 7,174.57 & 7,109.62 & 10,418.87 & \textbf{6,994.77} & 7,102.57 & 7,046.37 & 5,631.00 \\
CFLP\_10\_10\_500 & 7,171.79 & 7,068.91 & 10,410.19 & \textbf{7,003.30} & 7,102.57 & 7,084.46 & 5,643.68 \\
CFLP\_10\_10\_1000 & 7,154.60 & \textbf{7,040.70} & 10,406.08 & 7,088.56 & 7,102.57 & 7,064.36 & 5,622.40 \\
CFLP\_25\_25\_100 & \textbf{11,773.01} & \textbf{11,773.01} & 23,309.73 & 11,864.83 & 11,811.39 & 12,100.73 & 10,312.21 \\
CFLP\_25\_25\_500 & \textbf{11,726.34} & \textbf{11,726.34} & 23,310.34 & 12,170.67 & 11,811.39 & 12,051.51 & 10,277.01 \\
CFLP\_25\_25\_1000 & \textbf{11,709.90} & \textbf{11,709.90} & 23,309.85 & 11,868.04 & 11,811.39 & 12,041.12 & 10,263.37 \\
CFLP\_50\_50\_100 & 25,236.33 & \textbf{25,019.64} & 45,788.45 & 25,349.21 & 26,309.43 & 26,004.88 & 18,290.63 \\
CFLP\_50\_50\_500 & 25,281.13 & \textbf{24,964.33} & 45,786.97 & 28,037.66 & 26,287.48 & 25,986.50 & 18,209.77 \\
CFLP\_50\_50\_1000 & 25,247.77 & \textbf{24,981.70} & 45,787.18 & 30,282.41 & 26,309.43 & 26,002.78 & 18,217.14 \\
\bottomrule
\end{tabular}}
\caption{CFLP detailed objective results: each row represents an average over 10 2SP instance with varying scenario sets.  
``True objective" for \{\EVA{},\PSA{},LR\} is the cost of the first-stage solution obtained from the approximate MIP evaluated on the second-stage scenarios.  For EF it is the objective at the solving limit.  
``Approximate-MIP objective" is objective from the approximate MIP for \{\EVA{},\PSA{},LR\}. All times in seconds. 
}
\label{tab:obj_CFLP}
\end{table*}

\begin{table*}[t]\centering\resizebox{0.78\textwidth}{!}{
\begin{tabular}{l|rrrr|rrr}
\toprule
Problem & \multicolumn{4}{c}{True objective} & \multicolumn{3}{c}{Approximate-MIP objective} \\
\cmidrule(lr){2-5}
\cmidrule(lr){6-8}
{} & {\EVA{}} & {\PSA{}} & {LR} & {EF} & {\EVA{}\ } & {\PSA{}\ } & {LR\ } \\
\midrule
SSLP\_10\_50\_50 & -354.96 & -354.96 & -63.00 & \textbf{-354.96} & -350.96 & -339.42 & -294.69 \\
SSLP\_10\_50\_100 & \textbf{-345.86} & \textbf{-345.86} & -49.62 & -345.86 & -350.96 & -328.54 & -283.96 \\
SSLP\_10\_50\_500 & \textbf{-349.54} & \textbf{-349.54} & -54.68 & -349.54 & -350.96 & -332.82 & -288.02 \\
SSLP\_10\_50\_1000 & \textbf{-350.07} & \textbf{-350.07} & -55.45 & -235.22 & -350.96 & -333.46 & -288.55 \\
SSLP\_10\_50\_2000 & \textbf{-350.07} & \textbf{-350.07} & -54.72 & -172.73 & -350.96 & -332.87 & -288.19 \\
SSLP\_15\_45\_5 & -247.27 & -206.83 & -249.51 & \textbf{-255.55} & -238.44 & -259.11 & -58.28 \\
SSLP\_15\_45\_10 & -249.58 & -209.49 & -252.89 & \textbf{-257.41} & -238.44 & -265.92 & -64.01 \\
SSLP\_15\_45\_15 & -251.10 & -208.86 & -254.58 & \textbf{-257.68} & -238.44 & -267.01 & -66.71 \\
SSLP\_5\_25\_50 & -125.22 & -123.15 & 14.50 & \textbf{-125.36} & -121.64 & -110.18 & -119.98 \\
SSLP\_5\_25\_100 & -120.91 & -119.03 & 19.87 & \textbf{-120.94} & -121.64 & -109.59 & -117.79 \\
\bottomrule
\end{tabular}}
\caption{SSLP detailed objective results: each row represents an average over eleven 2SP instance with varying scenario sets.  See Table~\ref{tab:obj_CFLP} for a detailed description of the columns.}
\label{tab:obj_SSLP}
\end{table*}

\begin{table*}[t]\centering\resizebox{0.78\textwidth}{!}{
\begin{tabular}{l|rrrr|rrr}
\toprule
Problem & \multicolumn{4}{c}{True objective} & \multicolumn{3}{c}{Approximate-MIP objective} \\
\cmidrule(lr){2-5}
\cmidrule(lr){6-8}
{} & {\EVA{}} & {\PSA{}} & {LR} & {EF} & {\EVA{}\ } & {\PSA{}\ } & {LR\ } \\
\midrule
INVP\_B\_E\_4 & -51.56 & -55.29 & -46.25 & \textbf{-57.00} & -58.59 & -52.15 & -63.67 \\
INVP\_B\_E\_9 & -54.86 & -58.15 & -53.11 & \textbf{-59.33} & -58.81 & -55.33 & -63.67 \\
INVP\_B\_E\_36 & -59.55 & -58.19 & -58.86 & \textbf{-61.22} & -59.38 & -57.92 & -63.67 \\
INVP\_B\_E\_121 & -61.44 & -60.78 & -61.06 & \textbf{-62.29} & -59.60 & -58.91 & -63.67 \\
INVP\_B\_E\_441 & -59.60 & -59.83 & -59.91 & \textbf{-61.32} & -59.91 & -58.51 & -63.67 \\
INVP\_B\_E\_1681 & -59.81 & - & -59.30 & \textbf{-60.63} & -59.94 & - & -63.67 \\
INVP\_B\_E\_10000 & \textbf{-59.85} & - & -58.68 & -58.98 & -59.95 & - & -63.67 \\
INVP\_B\_H\_4 & -51.75 & -51.36 & -51.75 & \textbf{-56.75} & -58.12 & -52.41 & -61.24 \\
INVP\_B\_H\_9 & -56.56 & -56.56 & -56.56 & \textbf{-59.56} & -61.78 & -56.67 & -61.24 \\
INVP\_B\_H\_36 & -59.31 & -59.31 & -59.31 & \textbf{-60.28} & -59.38 & -59.52 & -61.24 \\
INVP\_B\_H\_121 & -59.93 & -59.93 & -59.93 & \textbf{-61.01} & -60.22 & -60.54 & -61.24 \\
INVP\_B\_H\_441 & -60.14 & -58.07 & -60.14 & \textbf{-61.44} & -60.23 & -58.13 & -61.24 \\
INVP\_B\_H\_1681 & \textbf{-60.47} & - & \textbf{-60.47} & -60.04 & -60.57 & - & -61.24 \\
INVP\_B\_H\_10000 & \textbf{-60.53} & - & \textbf{-60.53} & -58.93 & -60.65 & - & -61.24 \\
INVP\_I\_E\_4 & -55.35 & \textbf{-63.50} & -52.50 & \textbf{-63.50} & -66.79 & -58.96 & -71.57 \\
INVP\_I\_E\_9 & -61.63 & -64.80 & -61.89 & \textbf{-66.56} & -66.70 & -61.70 & -71.57 \\
INVP\_I\_E\_36 & -66.03 & -66.25 & -67.08 & \textbf{-69.86} & -67.39 & -65.18 & -71.57 \\
INVP\_I\_E\_121 & -67.35 & -67.92 & -69.07 & \textbf{-71.12} & -67.39 & -66.70 & -71.57 \\
INVP\_I\_E\_441 & -67.55 & -69.16 & -67.39 & \textbf{-69.64} & -67.63 & -67.43 & -71.57 \\
INVP\_I\_E\_1681 & -67.95 & -66.73 & -66.52 & \textbf{-68.85} & -67.69 & -67.62 & -71.57 \\
INVP\_I\_E\_10000 & \textbf{-67.94} & - & -65.67 & -67.04 & -67.82 & - & -71.57 \\
INVP\_I\_H\_4 & -54.75 & -55.78 & -54.75 & \textbf{-63.50} & -65.31 & -59.99 & -66.07 \\
INVP\_I\_H\_9 & -59.78 & -65.25 & -59.78 & \textbf{-65.78} & -64.15 & -61.08 & -66.07 \\
INVP\_I\_H\_36 & -63.78 & -64.80 & -63.78 & \textbf{-67.11} & -66.79 & -63.76 & -66.07 \\
INVP\_I\_H\_121 & -65.03 & -64.37 & -65.03 & \textbf{-67.75} & -65.38 & -64.64 & -66.07 \\
INVP\_I\_H\_441 & -65.12 & -65.12 & -65.12 & \textbf{-67.24} & -67.13 & -65.16 & -66.07 \\
INVP\_I\_H\_1681 & \textbf{-65.63} & -65.34 & \textbf{-65.63} & -65.41 & -65.87 & -65.03 & -66.07 \\
INVP\_I\_H\_10000 & \textbf{-65.66} & - & \textbf{-65.66} & -64.63 & -66.45 & - & -66.07 \\
\bottomrule
\end{tabular}}
\caption{INVP detailed objective results: each row represents single instance.  See Table~\ref{tab:obj_CFLP} for a detailed description of the columns.}
\label{tab:obj_INVP}
\end{table*}

\begin{table*}[t]\centering\resizebox{0.7\textwidth}{!}{
\begin{tabular}{l|rrrr|rrr}
\toprule
Problem & \multicolumn{4}{c}{True objective} & \multicolumn{3}{c}{Approximate-MIP objective} \\
\cmidrule(lr){2-5}
\cmidrule(lr){6-8}
{} & {\EVA{}} & {\PSA{}} & {LR} & {EF} & {\EVA{}\ } & {\PSA{}\ } & {LR\ } \\
\midrule
PP\_125 & 264.30 & 173.10 & -10.00 & \textbf{273.19} & 189.75 & 177.12 & 70.75 \\
PP\_216 & 200.29 & 131.83 & -10.00 & \textbf{220.25} & 189.75 & 168.10 & 70.75 \\
PP\_343 & 206.38 & 122.90 & -10.00 & \textbf{207.77} & 189.75 & 172.17 & 70.75 \\
PP\_512 & 204.41 & 134.83 & -10.00 & \textbf{223.86} & 189.75 & 162.54 & 70.75 \\
PP\_729 & 219.42 & 137.97 & -10.00 & \textbf{222.48} & 189.75 & 167.55 & 70.75 \\
PP\_1000 & 202.50 & 126.30 & -10.00 & \textbf{215.25} & 189.75 & 165.59 & 70.75 \\
\bottomrule
\end{tabular}}
\caption{PP detailed objective results: each row represents single instance.  See Table~\ref{tab:obj_CFLP} for a detailed description of the columns.}
\label{tab:obj_PP}
\end{table*}

\section{SSLP SIPLib Results}
\label{app:siplib}

In this section we report optimally gaps and solving times on the publicly available SSLP SIPLib instances in Table~\ref{tab:res_siplib}.  From the table, we can see that both \EVA{} and \PSA{} do quite well in terms of finding solutions, especially in the larger scenario case where they obtain optimal first-stage solutions.  
Perhaps, the most impressive results here is that \EVA{} is able to obtain optimal results for many instances in $\sim$0.1 seconds.

\begin{table*}[t]\centering\resizebox{0.68\textwidth}{!}{
\begin{tabular}{l|rrr|rrr}
\toprule
Problem & \multicolumn{3}{c}{Gap to Optimal (\%)} & \multicolumn{3}{c}{Solving Time} \\
\cmidrule(lr){2-4}
\cmidrule(lr){5-7}
{} & {\EVA{}} & {\PSA{}} & {EF} & {\EVA{}\ } & {\PSA{}\ } & {EF\ } \\
\midrule
SSLP\_10\_50\_50 & \textbf{0.00} & \textbf{0.00} & \textbf{0.00} & \textbf{0.11} & 4.83 & 10,801.27 \\
SSLP\_10\_50\_100 & \textbf{0.00} & \textbf{0.00} & \textbf{0.00} & \textbf{0.11} & 11.66 & 10,800.04 \\
SSLP\_10\_50\_500 & \textbf{0.00} & \textbf{0.00} & 0.00 & \textbf{0.11} & 107.88 & 10,818.23 \\
SSLP\_10\_50\_1000 & \textbf{0.00} & \textbf{0.00} & 28.64 & \textbf{0.12} & 383.51 & 10,800.26 \\
SSLP\_10\_50\_2000 & \textbf{0.00} & \textbf{0.00} & 51.24 & \textbf{0.13} & 4,523.19 & 10,800.20 \\
SSLP\_15\_45\_5 & 0.46 & 19.59 & \textbf{0.00} & 0.32 & \textbf{0.28} & 4.17 \\
SSLP\_15\_45\_10 & 1.57 & 18.23 & \textbf{0.00} & \textbf{0.25} & 0.40 & 3.71 \\
SSLP\_15\_45\_15 & 0.53 & 16.51 & \textbf{0.00} & \textbf{0.41} & 0.72 & 4.74 \\
SSLP\_5\_25\_50 & \textbf{0.00} & 2.15 & \textbf{0.00} & \textbf{0.26} & 0.92 & 2.35 \\
SSLP\_5\_25\_100 & \textbf{0.00} & 1.40 & \textbf{0.00} & \textbf{0.18} & 1.68 & 8.87 \\
\bottomrule
\end{tabular}}
\caption{SSLP SIPLib gap and time comparison among methods. Optimal SIPLib instances values obtained from \cite{ahmed2015siplib}.  
``Gap to Optimal" is the percent gap to the optimal solution. 
``Solving Time" is the solving to of the approximate MIP and EF.
All times in seconds.}
\label{tab:res_siplib}
\end{table*}

\section{Model \& Dataset Selection}
\label{app:ms_ds}

\subsection{Model Selection}

For the supervised learning task, we implement linear regression using Scikit-learn 1.0.1 \citep{scikit-learn}.  In this case we use the base estimator with no regularization.  
The \PSA{}/\EVA{} neural models are all implemented using Pytorch 1.10.0 \citep{NEURIPS2019_9015}.  
For model selection, we use random search over 100 configurations for each problem setting.  For \PSA{} and \EVA{} we sample configurations from Table~\ref{tab:random_search_space}.  For both cases we limit the ReLU layers to a single layer with a varying hidden dimension.  In the \PSA{} case the choice of the ReLU hidden dimension is limited since a large number of predictions each with a large hidden dimension can lead to MILPs which are prohibitively expensive to solve.  For the \EVA{} specific hidden dimensions, we have 3 layers, with Embed hidden dimension 1 and Embed hidden dimension 2 corresponding to layers before the aggregation and Embed hidden dimension 3 being a final hidden layer after the aggregation.  

\begin{figure}[ht]
\begin{center}
\centerline{\includegraphics[scale=0.5]{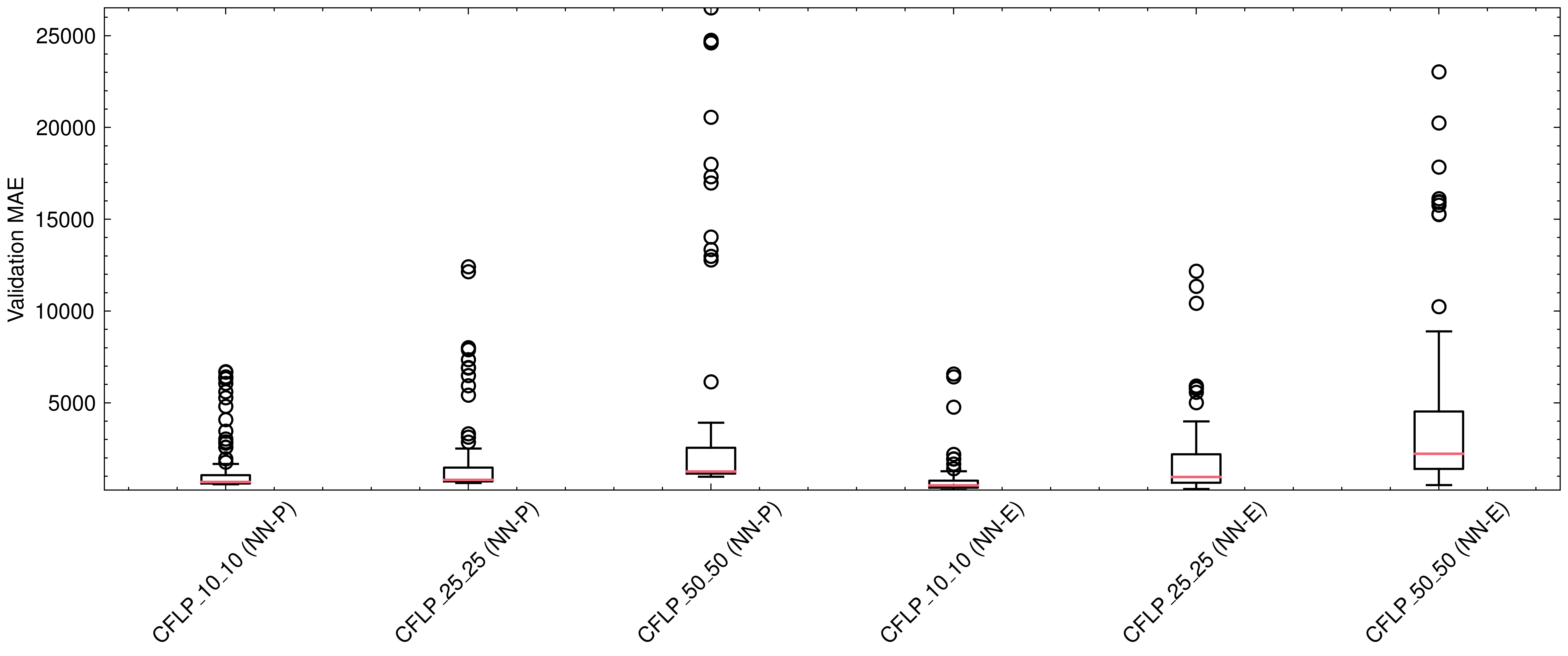}}
\caption{CFLP validation MAE over random search configurations for \PSA{} and \EVA{} models.   }
\label{fig:cflp_rs_box_plots}
\end{center}
\end{figure}

\begin{figure}[ht]
\begin{center}
\centerline{\includegraphics[scale=0.5]{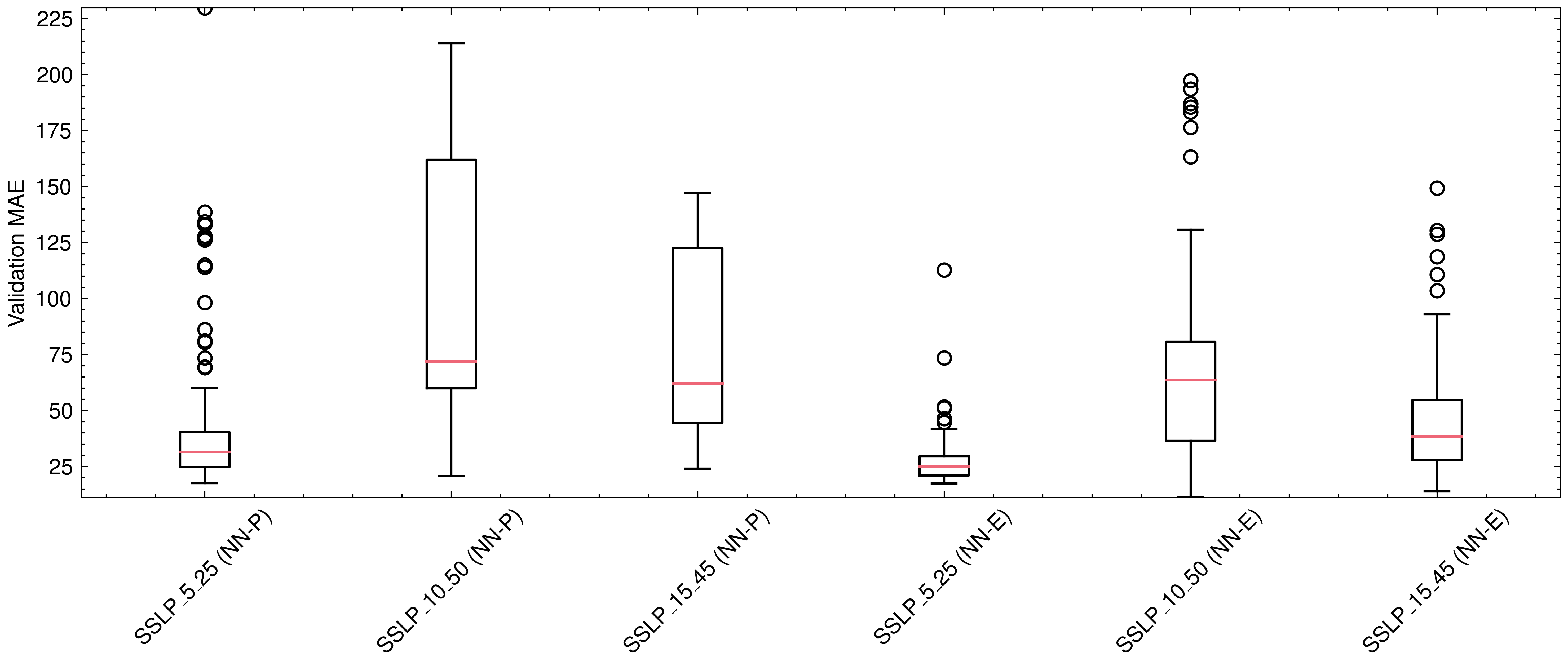}}
\caption{SSLP validation MAE over random search configurations for \PSA{} and \EVA{} models.  }
\label{fig:sslp_rs_box_plots}
\end{center}
\end{figure}

\begin{figure}[ht]
\begin{center}
\centerline{\includegraphics[scale=0.5]{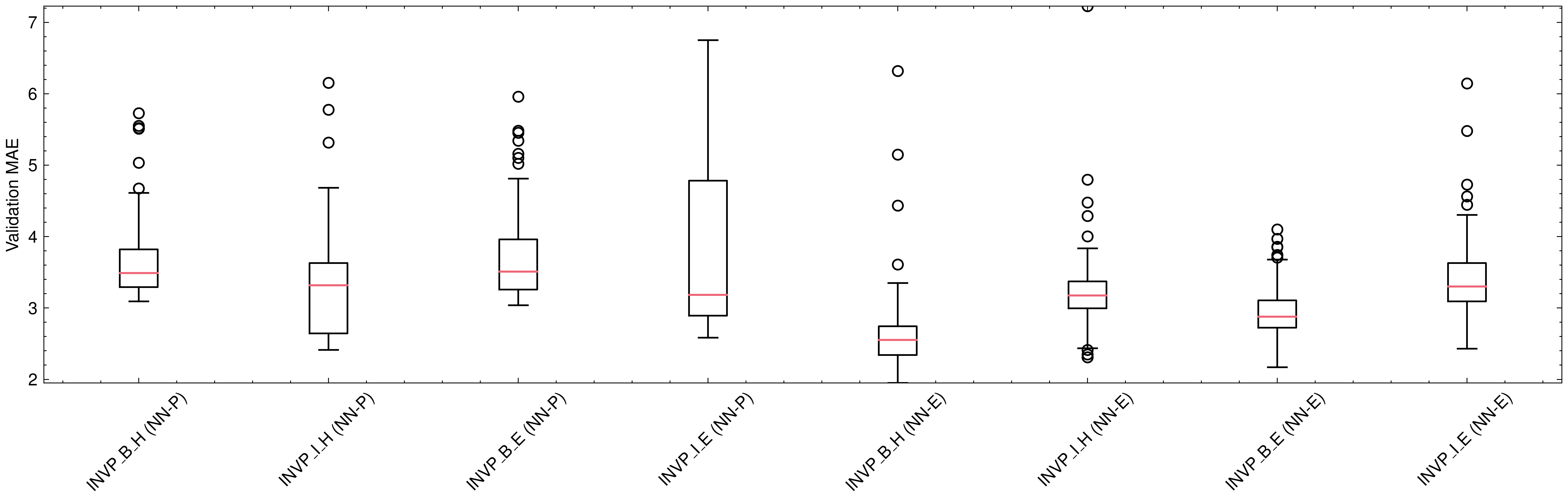}}
\caption{INVP validation MAE over random search configurations for \PSA{} and \EVA{} models.  }
\label{fig:ip_rs_box_plots}
\end{center}
\end{figure}

\begin{figure}[ht]
\begin{center}
\centerline{\includegraphics[scale=0.5]{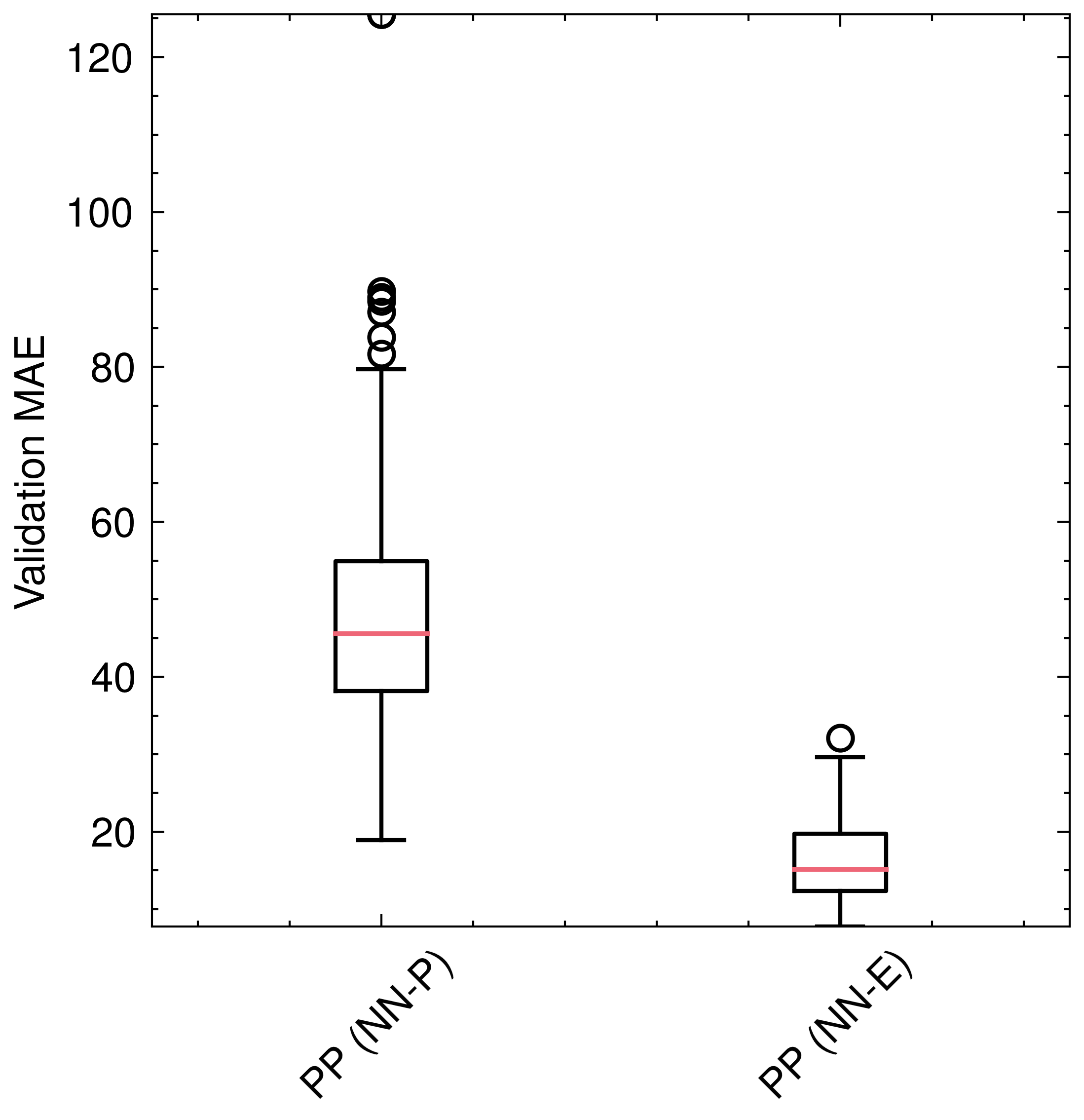}}
\caption{PP validation MAE over random search configurations for \PSA{} and \EVA{} models.  }
\label{fig:pp_rs_box_plots}
\end{center}
\end{figure}

In Tables~\ref{tab:mc_rs_best} and \ref{tab:sc_rs_best} we report the best parameters for each problem setting for the \PSA{} and \EVA{} models, respectively.  In addition, we report the validation MSE across all 100 configurations for each problem in box plots in Figures~\ref{fig:cflp_rs_box_plots} to \ref{fig:pp_rs_box_plots}.  From the box plots we can observe that lower validation MAE configurations are quite common as the medians are typically not too far from the lower tails of the distributions.  This indicates that hyperparameter selection can be helpful when attempting to improve the second-stage cost estimates, however, the gains are marginal in most cases.

\begin{table*}[t]\centering\resizebox{0.8\textwidth}{!}{
\begin{tabular}{l|cc}
\toprule
{Parameter}                     & {\PSA{}}                          & {\EVA{} } \\
\midrule
    Batch size                  & $\{16, 32, 64, 128\}$         & $\{16, 32, 64, 128\}$  \\
    Learning rate               & $[1e^{-5},1e^{-1}]$           & $[1e^{-5},1e^{-1}]$ \\
    L1 weight penalty           & $[1e^{-5}, 1e^{-1}]$          & $[1e^{-5}, 1e^{-1}]$  \\
    L2 weight penalty           & $[1e^{-5}, 1e^{-1}]$          & $[1e^{-5}, 1e^{-1}]$\\
    Optimizer                   & \{Adam, Adagrad, RMSprop\}    & \{Adam, Adagrad, RMSprop\}\\
    Dropout                     & $[0, 0.5]$                    & $[0, 0.5]$ \\     
    \# Epochs                   & $1000$                        & $2000$ \\
    ReLU hidden dimension       & $\{32, 64\}$                  & $\{64, 128, 256, 512\}$  \\
    Embed hidden dimension 1    & -                             & $\{64, 128, 256, 512\}$   \\
    Embed hidden dimension 2    & -                             & $\{16, 32, 64, 128\}$  \\
    Embed hidden dimension 3    & -                             & $\{8, 16, 32, 64\}$  \\
\bottomrule
\end{tabular}}
\caption{Random search parameter space for \PSA{} and \EVA{} models.  For values in \{\}, we sample with equal probability for each discrete choice.  For values in [], we sample a uniform distribution with the given bounds.  For single values, we keep it fixed across all configurations.  A value of - means that parameter is not applicable for the given model type. }
\label{tab:random_search_space}
\end{table*}


\begin{table*}[t]\centering\resizebox{\textwidth}{!}{
\begin{tabular}{l|rrrrrrrrrrr}
\toprule
{Parameter}                     & {CFLP\_10\_10}    & {CFLP\_25\_25}    & {CFLP\_50\_50}    & {SSLP\_5\_25}    & {SSLP\_10\_50}    & {SSLP\_15\_45} & {INVP\_B\_E} & {INVP\_B\_H} & {INVP\_I\_E} & {INVP\_I\_H} & {PP}\\
\midrule
    Batch size                  & 128       & 16        & 128       & 128       & 128     & 64 & 16 & 32 & 32 & 128 & 64\\   
    Learning rate               & 0.05029   & 0.05217   & 0.00441   & 0.03385   & 0.07015 &  0.08996 & 0.00435 & 0.00521 & 0.06613 & 0.01614 & 0.0032\\   
    L1 weight penalty           & 0.07512   & 0.00551   & 0.08945   & 0.03232   & 0.07079 &  0.09105 & 0.08321 & 0.05754 & 0.01683 & 0.
    01859 & 0 \\   
    L2 weight penalty           & 0.08343   & 0.02846   & 0.08602   & 0.0       & 0.01792 &  0.0 & 0.01047 & 0.02728 & 0 & 0& 0.0361\\   
    Optimizer                   & Adam      & Adam      & Adam      & RMSprop   & RMSprop &  RMSprop & RMSProp & RMSProp & Adam & Adam & Adam \\    
    Dropout                     & 0.02198   & 0.02259   & 0.05565   & 0.00914   & 0.01944 &  0.02257 & 0.17237 & 0.13698 & 0.04986 & 0.0859 & 0.05945\\   
    ReLU hidden dimension       & 64        & 32        & 64        & 32        & 64 &  32 & 64 & 64 & 64 & 32 & 64\\ 
\bottomrule
\end{tabular}}
\caption{\PSA{} best configurations from random search.}
\label{tab:mc_rs_best}
\end{table*}

\begin{table*}[t]\centering\resizebox{\textwidth}{!}{
\begin{tabular}{l|rrrrrrrrrrr}
\toprule
{Parameter}                     & {CFLP\_10\_10}    & {CFLP\_25\_25}    & {CFLP\_50\_50}    & {SSLP\_5\_25}    & {SSLP\_10\_50}    & {SSLP\_15\_45} & {INVP\_B\_E} & {INVP\_B\_H} & {INVP\_I\_E} & {INVP\_I\_H} & {PP}\\
\midrule
    Batch size                  & 32        & 16        & 128       & 64        & 64 & 32 & 128 & 32 & 16 & 128 & 64 \\   
    Learning rate               & 0.0263    & 0.06571   & 0.02906   & 0.08876   & 0.07633 & 0.03115  & 0.01959 & 0.00846 & 0.02841 & 0.02867  & 0.08039 \\   
    L1 weight penalty           & 0.02272   & 0.06841   & 0.05792   & 0.0       & 0.04529   & 0.07182 & 0.0 & 0.0 & 0.00022 & 0 & 0\\   
    L2 weight penalty           & 0.05747   & 0.0       & 0.04176   & 0.03488   & 0.0       & 0.0& 0 & 0.09007 & 0.02272 & 0.01882 & 0\\   
    Optimizer                   & RMSprop   & Adam      & Adam      & Adam      & RMSprop   & Adam & Adagrad & Adam & Adagrad & Adagrad & Adam\\    
    Dropout                     & 0.01686   & 0.0028    & 0.03318   & 0.00587   & 0.00018   & 0.0088 & 0.08692 & 0.04096 & 0.01854 & 0.01422 & 0.0072\\   
    ReLU hidden dimension       & 128       & 256       & 256       & 256       & 64        & 256 & 256 & 256 & 256 & 256 & 512\\  
    Embed hidden dimension 1    & 512       & 128       & 256       & 64        & 128       & 512 & 256 & 512 & 64 & 256 & 512\\ 
    Embed hidden dimension 2    & 16        & 64        & 64        & 16        & 32        & 64 & 16 & 16 & 32 & 32 & 128\\ 
    Embed hidden dimension 3    & 16        & 16        & 8         & 32        & 64        & 16 & 32 & 16 & 8 & 64 & 16\\ 
\bottomrule
\end{tabular}}
\caption{\EVA{} best configurations from random search.}
\label{tab:sc_rs_best}
\end{table*}

\subsection{Dataset Size Selection}

\begin{figure}[h!]
  \centering
  \begin{minipage}[b]{0.45\textwidth}
    \includegraphics[width=\textwidth]{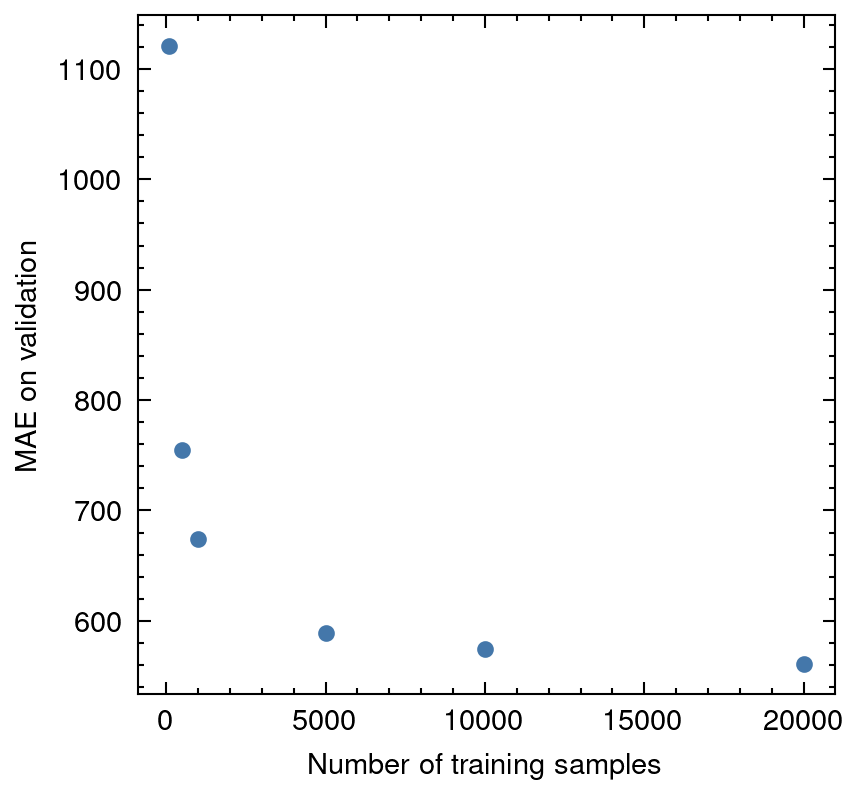}
    \caption{\PSA{} dataset sizing results}
    \label{fig:dss_mc}
  \end{minipage}
  \hfill
  \begin{minipage}[b]{0.45\textwidth}
    \includegraphics[width=\textwidth]{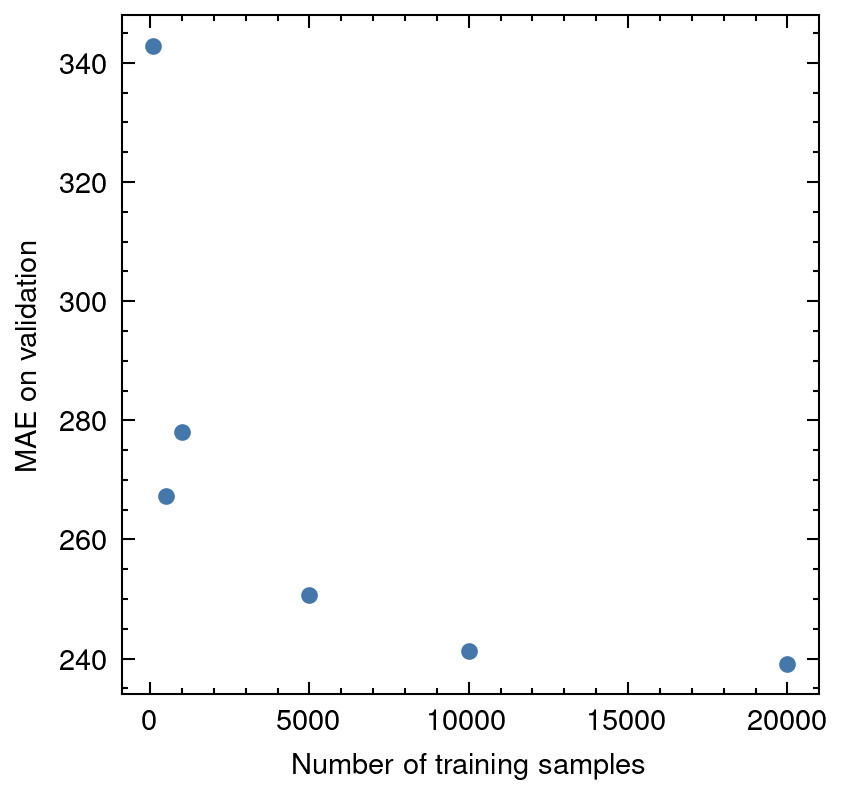}
    \caption{\EVA{} dataset sizing results}
    \label{fig:dss_sc}
  \end{minipage}
\end{figure}

In this section, we report results for varying dataset sizes.  Here, we report only the results for a single problem setting, namely, CFLP\_10\_10.  We use a validation with 5000 samples and training sets with 100, 500, 1000, 5000, 10000, and 20000 samples.  Model selection with random search is done for each training set size as described in the previous section.  Figures~\ref{fig:dss_mc} and \ref{fig:dss_sc} report the results for the \PSA{} and \EVA{} models respectively.  In both cases, we can see an improvement in validation MAE with increases in the dataset sizes, however, diminishing returns start to occur when increasing the number of samples above 5000 samples.  This motivates the choice of dataset sizes which we use in the remainder of the experiments.  Specifically, we use 10000 samples for the \PSA{} case as data generation is quite fast.  For the \EVA{} case we limit the number of samples to 5000 as we only see a small improvement of \%4 in validation MAE at the cost of doubling the compute time.

\end{document}